\documentclass[a4paper,leqno]{article}

\usepackage{amsmath}
\usepackage{amssymb}
\usepackage{amsfonts}
\usepackage{enumerate}
\usepackage{theorem}

\theoremstyle{plain}
\newtheorem{teo}{Theorem}[section]
\newtheorem{lem}[teo]{Lemma}
\newtheorem{cor}[teo]{Corollary}
\newtheorem{prop}[teo]{Proposition}

{\theorembodyfont{\rmfamily}
\newtheorem{defin}[teo]{Definition}
\newtheorem{oss}[teo]{Remark}
\newtheorem{exam}[teo]{Example}
}

\DeclareMathOperator*{\argmin}{arg\,min}

\renewcommand{\eqref}[1]{\textnormal{(\ref{#1})}}

\numberwithin{equation}{section}

\newcommand{\cvd}{\hfill$\square$}
\newcommand{\proof}[1]{\noindent\textsc{Proof#1}}
\newcommand{\rmi}{\mathrm{i}}

\newcommand{\rmd}{\mathrm{d}}

\title{Discrete approximation and regularisation for the inverse conductivity problem}

\author{Luca \textsc{Rondi}\\
\normalsize{Dipartimento di Matematica e Geoscienze}\\
\normalsize{Universit\`a degli Studi di Trieste}\\
\normalsize{via Valerio, 12/1}\\
\normalsize{34127 Trieste, Italy}\\
\normalsize{\texttt{rondi@units.it}}}

\date{
}

\begin{document}

\maketitle

\setcounter{section}{0}
\setcounter{secnumdepth}{2}

\begin{abstract}
We study the inverse conductivity problem with discontinuous conductivities.
We consider, simultaneously, a regularisation and a discretisation for a variational approach to solve the inverse problem. We show that, under suitable choices of the regularisation and discretisation parameters, the discrete regularised solutions converge, as the noise level on the measurements goes to zero, to the looked for solution of the inverse problem.

\medskip


\medskip

\end{abstract}

\bigskip

\begin{flushright}
\emph{Dedicated to Giovanni Alessandrini on the occasion of his 60th birthday}
\end{flushright}

\section{Introduction}

In this paper we consider the inverse conductivity problem with discontinuous conductivity. For a given conducting body contained in a bounded domain $\Omega\subset \mathbb{R}^N$, $N\geq 2$,
we call $X$ the space of admissible conductivities, or better conductivity tensors, in $\Omega$. For any $\sigma\in X$, we call $\Lambda(\sigma)$ either the Dirichlet-to-Neumann map, or the Neumann-to-Dirichlet map, corresponding to $\sigma$. It is a well-known fact that $\Lambda(\sigma)$ is a bounded linear operator between suitable
Banach spaces defined on the boundary of $\Omega$, and we call $Y$ the space of these bounded linear operators. The forward operator $\Lambda:X\to Y$ is the one that to each $\sigma\in X$ associates $\Lambda(\sigma)\in Y$.

The aim of the inverse problem is to determine an unknown conductivity in $\Omega$ by performing suitable electrostatic measurements of current and voltage type on the boundary.
If $\sigma_0$ is the conductivity we aim to recover by solving our inverse problem, then 
we measure its corresponding $\Lambda(\sigma_0)\in Y$. Due to the noise that is present in the measurements, actually the information that we are able to collect is
$\hat{\Lambda}\in Y$, which is a perturbation of $\Lambda(\sigma_0)$. We call
$\|\hat{\Lambda}-\Lambda(\sigma_0)\|_Y$ the noise level of the measurements and we notice that the choice of the space $Y$ corresponds to the way we measure the errors in our measurements.

The inverse problem may be stated, at least formally, in the following way.
Given our measurements $\hat{\Lambda}$, we wish to find $\sigma\in X$ such that
\begin{equation}\label{formalinvpbm}
\Lambda(\sigma)=\hat{\Lambda}.
\end{equation}
Due to the noise, such a problem may not have any solution, therefore we better consider a least-square formulation
\begin{equation}\label{formalinvpbmbis}
\min_{\sigma\in X}\|\Lambda(\sigma)-\hat{\Lambda}\|_Y.
\end{equation}
Unfortunately, the inverse conductivity problem is ill-posed, therefore to solve 
\eqref{formalinvpbmbis} numerically, a regularisation strategy need to be implemented.
Considering a regularisation \`a la Tikhonov, this means to choose a regularisation operator $R$, usually a norm or a seminorm, and a regularisation parameter $a$ and solve
\begin{equation}\label{formalinvpbmreg}
\min_{\sigma\in X}\|\Lambda(\sigma)-\hat{\Lambda}\|_Y+aR(\sigma) .
\end{equation}
A solution to \eqref{formalinvpbmreg} is called a regularised solution. A good regularisation operator need to satisfy the following two criteria. First of all, it should make the minimisation process stable from a numerical point of view. Second, the regularised solution should be a good approximation of the looked for solution of the inverse problem.

For the nonsmooth case, often this second requirement is not proved analytically but rather it is (not rigorously) justified by
numerical tests only. However, a convergence analysis, using techniques inspired by variational convergences such as
$\Gamma$-convergence, allows to rigorously justify the choice of the regularisation operator, \cite{Ron08}. For the inverse conductivity problem with discontinuous conductivity, by this technique, in the same paper \cite{Ron08}, the use of some of the usually employed regularisation methods  was rigorously justified. For instance, a convergence analysis was developed for regularisations such as the total variation penalisation or the Mumford-Shah functional. Several other works followed this approach, for instance it was extended to smoothness or sparsity penalty regularisations for the inverse conductivity problem in \cite{Ji-Ma}, whereas in \cite{Jia-Ma-Pa} the analysis for the Mumford-Shah functional was slightly refined and applied to other inverse problems.

Once the regularisation operator is chosen, and proved to be effective, the issue of the numerical approximation for the regularised problem comes into play. One of the key points of the numerical approximation is represented by the discretisation of the regularised minimum problem. Again, two issues come forward. The first one is the choice of the kind of space of discrete unknowns we intend to use. The second 
important issue is how fine the discretization should be. A compromise is necessary between a better resolution
(finer discretization) and a more stable reconstruction (coarser discretization). Again, the discrete regularised solution, that is, the solution to the regularised problem \eqref{formalinvpbmreg} with $\sigma$ varying in such a discrete subset, should
be a good approximation of the solution of the inverse problem. Actually, for inverse problems,  this may not be necessarily so, as an example in
\cite{Riv-Bar-Ob} shows. Therefore, studying the
effect of the discretisation when solving an inverse problem is not at all an easy task.
This fundamental and nontrivial issue went rather unlooked, at least for the inverse conductivity problem and other classical inverse problems dealing with nonsmooth unknowns.

The crucial point we wish to address here is the following. We want to simultaneously fix both the regularisation parameter and the discretisation parameter, in correspondence to the given noise level, such that the discrete regularised solutions converge, as the noise level goes to zero, to the solution of the inverse problem.
Previously, only the analysis of the approximation of the
regularised problem with discrete ones, with a fixed regularisation parameter, was performed. For instance, a nice finite element approximation for the inverse conductivity problem, with the total variation as regularisation, may be found in \cite{Ge-Ji-Lu}.
In \cite{Ron-San}, instead, it was proved that the regularised inverse conductivity problem, with the Mumford-Shah as a regularisation term, could be well approximated by replacing the Mumford-Shah with its approximating Ambrosio-Tortorelli functionals developed in \cite{Amb-Tor90,Amb-Tor92}. Here the approximating parameter for the Ambrosio-Tortorelli functionals may be seen as another version of the discretisation parameter.

Actually, the first attempt to vary, in a suitable way, the
regularisation and discretisation parameters simultaneously, may be found in   
a Master thesis supervised by the author, \cite{Che}. There the Ambrosio-Tortorelli functionals were considered, and their approximating parameter and the regularisation parameter were chosen accordingly to the noise level to guarantee the required convergence of this type of regularised solutions. For the convenience of the reader, we present a brief summary of this result in Subsection~\ref{cherinisubs} of the present paper.

The main result of the paper, Theorems~\ref{mainthm} and \ref{2dimcasebis}, is contained in Subsection~\ref{discretesubsec}. We consider the inverse conductivity problem and its regularisation by a total variation penalisation. We consider a discrete subset of admissible conductivities which is simply given by standard conforming piecewise linear finite elements over a regular triangulation. The triangulation is characterised by a discretisation parameter $h$, which is an upper bound for the diameter of any simplex forming the triangulation.

We show that, if we choose the regularisation parameter $a$ and the discretisation parameter $h$ according to the noise level, then the discrete regularised solution  
would converge to a solution of the inverse problem. An interesting feature of this result is that it shows that the discretisation parameter should go to zero in a polynomial way with respect to the noise level.

We remark that in this paper we limit ourselves to a very simple scenario but we believe that this is just a first step to tackle a full discretisation of the inverse conductivity problem, in a more general setting as well. This will be the object of future work.

It would also be very interesting to address the issue of convergence estimates.
In the smooth case they may be obtained by using Tikhonov regularisation for nonlinear operators, see for instance \cite{Eng-et-al}. Actually, for the inverse conductivity problem in the smooth case, some convergence estimates are available for the regularised solutions, without adding the discretisation, see for instance \cite{LMP} and \cite{Ji-Ma}. We notice that our technique involves $\Gamma$-convergence, which is of a qualitative nature thus does not lead easily to convergence estimates.

Finally we wish to mention that, for discrete sets of unknowns, that is, unknowns depending on a finite number of parameters, 
the usual ill-posedness of these kinds of inverse problems considerably reduces.
In fact, Lipschitz stability estimates may be obtained instead of the classical logarithmic
ones. Such an important line of research was initiated in \cite{Ale-Ves} and pursued in several other paper (let us mention the recent one \cite{A-dH-G-Sin} which is the closest to the setting we use in this paper). Unfortunately, the behaviour of the Lipschitz constant as the discretisation parameter approaches zero is extremely bad, as it explodes exponentially with respect to $h$, a fact firstly noted in \cite{Ron06}. This fact seems to prevent the use of these kinds of estimates at the discrete level to prove convergence estimate, or even just convergence, of discrete regularised solutions.

The plan of the paper is the following.
In Section~\ref{sec2}, besides fixing the notation and stating the inverse conductivity problem, we present a rather complete introduction to the regularisation issue for this inverse problem. Most of the material here is not new, a part from a few instances that we point out in a while, but our aim is to present a self-contained review to this line of research that is scattered in several papers. We begin with uniqueness results for scalar conductivities, that is, for the isotropic case, and nonuniqueness for symmetric conductivity tensors, that is, for the anisotropic case,
Subsection~\ref{statesubs}. We recall that nonuniqueness is due to the invariance of the boundary operators by smooth changes of variables of the domain $\Omega$ that keep fixed the boundary.

In Subsection~\ref{illsec}, we study the existence of a solution to
\eqref{formalinvpbm}. This part is mostly from \cite{Ron15}. We show that existence is
true in the anisotropic case, whereas it may fail in the isotropic case, see 
Example~\ref{noexistenceexam}. We notice that 
Example~\ref{noexistenceexam} appeared in a Master thesis supervised by the author, \cite{Des},
and it is a slight generalisation of a similar example in \cite{Ron15}. The crucial ingredient for both is a nice construction due to Giovanni that may be found in \cite[Example~4.4]{Ron15}.
Even if existence of \eqref{formalinvpbm} is guaranteed, the ill-posedness nature of this inverse problem implies that minimiser to \eqref{formalinvpbm} may fail to converge to the looked for solution to the inverse problem, as the noise level goes to zero. This is shown in three different examples,
Examples~\ref{examanisot}, \ref{AleCabexam} and \ref{examanisot2}.
Example~\ref{examanisot} shows how nonuniqueness in the anisotropic case leads to instability, see also Proposition~\ref{weakconprop} which is taken from \cite{Far-Kur-Rui} for a corresponding partial stability result.
Examples~\ref{AleCabexam} and \ref{examanisot2} deal with the isotropic case. The 
latter is new and slightly improves the former, which is taken from \cite{Ale-Cab}.

In Subsection~\ref{regularsubs}, we recall the approach to regularisation for inverse problems with nonsmooth unknowns, and in particular for the inverse conductivity problem with discontinuous conductivities, that was developed in \cite{Ron08}.

Section~\ref{discretesec} is the main of the paper. We investigate simultaneous numerical approximation and regularisation for the inverse conductivity problem with discontinuous conductivities. In 
Subsection~\ref{discretesubsec}, we present our main result, the convergence analysis of the discretisation by the finite element method coupled with a total variation regularisation. Finally,
in
Subsection~\ref{cherinisubs}, we present the result of \cite{Che}, that is, the convergence analysis for the regularisation by Ambrosio-Tortorelli functionals.

\subsubsection*{Acknowledgement}
The author is partially supported by Universit\`a degli Studi di Trieste through FRA 2014 and by GNAMPA, INdAM.

\section{Statement of the inverse problem, preliminary considerations, and previous results}\label{sec2}

Throughout the paper we shall keep fixed 
positive constants $\lambda_0$, $\lambda_1$, and $\tilde{\lambda}_1$, with $0<\lambda_0\leq  \lambda_1,\tilde{\lambda}_1$. The integer $N\geq 2$ will always denote the space dimension and we recall that we shall usually drop the dependence of any constant on $N$. For any Borel set $E\subset\mathbb{R}^N$, we denote with
$|E|$  its Lebesgue measure, whereas $\mathcal{H}^{N-1}(E)$ denotes its $(N-1)$-dimensional Hausdorff measure.

Throughout the paper we also fix $\Omega$, a bounded connected open set contained in $\mathbb{R}^N$, $N\geq 2$. We assume that $\Omega$ has a Lipschitz boundary in the following usual sense. For any $x\in\partial\Omega$ there exist $r>0$ and a Lipschitz function
$\varphi:\mathbb{R}^{N-1}\to \mathbb{R}$ such that, up to a rigid change of coordinates, we have
$$\Omega\cap B_r(x)=\{y=(y_1,\ldots,y_{N-1},y_N)\in B_r(x):\ y_N<\varphi(y_1,\ldots,y_{N-1})\}.$$

We call $\mathbb{M}^{N\times N}(\mathbb{R})$ the space of real valued $N\times N$ matrices. For any $\sigma\in \mathbb{M}^{N\times N}(\mathbb{R})$, with $N\geq 2$, several equivalent ellipticity conditions may be used. For example
\begin{equation}\label{ell2}
\left\{\begin{array}{ll}
\sigma\xi\cdot\xi\geq \lambda_0\|\xi\|^2&\text{ for any }\xi\in\mathbb{R}^N\\
\sigma^{-1}\xi\cdot\xi\geq \lambda_1^{-1}\|\xi\|^2&\text{ for any }\xi\in\mathbb{R}^N.\end{array}\right.
\end{equation}
Otherwise we can use
\begin{equation}\label{ell1}
\left\{\begin{array}{ll}
\sigma\xi\cdot\xi\geq \lambda_0\|\xi\|^2&\text{ for any }\xi\in\mathbb{R}^N\\
\|\sigma\|\leq \tilde{\lambda}_1
\end{array}\right.
\end{equation}
where $\|\sigma\|$ denotes its norm as a linear operator of $\mathbb{R}^N$ into itself.

The following remark shows that these two conditions are equivalent. If $\sigma$ satisfies \eqref{ell2} with constants $\lambda_0$ and $\lambda_1$, then it also satisfies \eqref{ell1} with constants $\lambda_0$ and $\tilde{\lambda}_1=\lambda_1$. If $\sigma$ satisfies \eqref{ell1} with constants $\lambda_0$ and $\tilde{\lambda}_1$, then it also satisfies \eqref{ell2} with constants $\lambda_0$ and $\lambda_1=\tilde{\lambda}_1^2/\lambda_0$. If $\sigma$ is symmetric then, picking $\tilde{\lambda}_1=\lambda_1$,
\eqref{ell2} and \eqref{ell1} are exactly equivalent and coincide with the condition
$$\lambda_0\|\xi\|^2\leq \sigma\xi\cdot\xi\leq \lambda_1\|\xi\|^2\quad \text{ for any }\xi\in\mathbb{R}^N,$$
that we write in short as follows
$$\lambda_0I_N\leq \sigma\leq \lambda_1I_N,$$
where $I_N$ is the $N\times N$ identity matrix.
Finally, if $\sigma=s I_N$, where $s$ is a real number, the condition simply reduces to
$$\lambda_0\leq s\leq \lambda_1.$$

We use the following classes of conductivity tensors in $\Omega$.
For positive constants $\lambda_0\leq\lambda_1$ we call
$\mathcal{M}(\lambda_0,\lambda_1)$
the set of $\sigma=\sigma(x)$, $x\in \Omega$, an $N\times N$ matrix whose entries are real valued  measurable functions in $\Omega$, such that, for almost any $x\in\Omega$, $\sigma(x)$ satisfies \eqref{ell2}. We call $\mathcal{M}_{sym}(\lambda_0,\lambda_1)$, respectively $\mathcal{M}_{scal}(\lambda_0,\lambda_1)$, the set of $\sigma\in \mathcal{M}(\lambda_0,\lambda_1)$ such that, for almost any $x\in\Omega$, $\sigma(x)$ is symmetric, respectively  $\sigma(x)=s(x) I_N$ with $s(x)$ a real number.
We say that $\sigma$ is a \emph{conductivity tensor} in $\Omega$
if $\sigma\in \mathcal{M}(\lambda_0,\lambda_1)$ for some constants $0<\lambda_0\leq\lambda_1$. We call $\mathcal{M}$ the class of conductivity tensors in $\Omega$.
We say that $\sigma$ is a \emph{symmetric conductivity tensor} in $\Omega$ if
$\sigma\in \mathcal{M}$ and $\sigma(x)$ is symmetric for almost any $x\in\Omega$.
We call $\mathcal{M}_{sym}$ the class of conductivity tensors in $\Omega$.
We say that $\sigma$ is a \emph{scalar conductivity} in $\Omega$ if
$\sigma\in \mathcal{M}$ and $\sigma(x)=s(x)I_N$, with $s(x)\in \mathbb{R}$, for almost any $x\in\Omega$.
We call $\mathcal{M}_{scal}$ the class of scalar conductivities in $\Omega$.

Since $\mathcal{M}\subset L^{\infty}(\Omega,\mathbb{M}^{N\times N}(\mathbb{R}))$,
we may measure the distance between any two conductivity tensors $\sigma_1$ and $\sigma_2$ in $\Omega$ with an $L^p$ metric, for any $p$, $1\leq p\leq +\infty$,
as follows
$$\|\sigma_1-\sigma_2\|_{L^p(\Omega)}=\|(\|\sigma_1-\sigma_2\|)\|_{L^p(\Omega)}.$$
With any of these $L^p$ metrics, any of the classes $\mathcal{M}(\lambda_0,\lambda_1),$
$\mathcal{M}_{sym}(\lambda_0,\lambda_1)$, and $\mathcal{M}_{scal}(\lambda_0,\lambda_1)$
is a complete metric space.

For any $p$, $1\leq p\leq +\infty$, we denote with $p'$ its conjugate exponent, that is $1/p+1/p'=1$.
For any $p$, $1< p< +\infty$,
we call
$W^{1-1/p,p}(\partial \Omega)$ the space of traces of
$W^{1,p}(\Omega)$ functions on $\partial \Omega$.
We recall that
$W^{1-1/p,p}(\partial \Omega)\subset L^p(\partial\Omega)$, with compact immersion.
For simplicity, we denote $H^1(\Omega)=W^{1,2}(\Omega)$,
$H^{1/2}(\partial\Omega)=W^{1/2,2}(\partial \Omega)$ and 
$H^{-1/2}(\partial\Omega)$ its dual.

We call $L^2_{\ast}(\partial\Omega)$ the subspace of functions $f\in L^2(\partial\Omega)$ such that $\int_{\partial\Omega}f=0$.
We set
$H^{-1/2}_{\ast}(\partial\Omega)$ the subspace of $g\in H^{-1/2}(\partial\Omega)$ such that
$$\langle g,1\rangle_{(H^{-1/2}(\partial \Omega),H^{1/2}(\partial \Omega))}=0.$$
We recall that $L^2_{\ast}(\partial\Omega)\subset H^{-1/2}_{\ast}(\partial\Omega)$, with compact immersion, if for any $g\in L^2_{\ast}(\partial\Omega)$ and 
any $\psi\in H^{1/2}(\partial \Omega)$ we define
\begin{equation}
\langle g,\psi\rangle_{(H^{-1/2}(\partial \Omega),H^{1/2}(\partial \Omega))}=
\int_{\partial\Omega}g\psi.
\end{equation}
Analogously, $H^{1/2}_{\ast}(\partial \Omega)$ is the subspace of $\psi\in H^{1/2}(\partial\Omega)$ such that
$\int_{\partial\Omega}\psi=0$. We have $H^{1/2}_{\ast}(\partial \Omega)\subset L^2_{\ast}(\partial\Omega)$, with compact immersion.
 
For any two Banach spaces $B_1$, $B_2$, $\mathcal{L}(B_1,B_2)$
will denote the Banach space of bounded linear operators from $B_1$ to $B_2$ with the usual operator norm.

\subsection{Statement of the problem and uniqueness results}\label{statesubs}

For any conductivity tensor $\sigma$ in $\Omega$, we define its Dirichlet-to-Neumann map
$DN(\sigma):H^{1/2}(\partial\Omega)\to H^{-1/2}(\partial\Omega)$
where for each $\varphi\in H^{1/2}(\partial\Omega)$,
$$DN(\sigma)(\varphi)[\psi]=\int_{\Omega}\sigma\nabla u\cdot\nabla \tilde{\psi} \quad\text{for any }\psi\in H^{1/2}(\partial\Omega)$$
where $u$ solves
\begin{equation}\label{Dirichletproblem1}
\left\{\begin{array}{ll}
\mathrm{div}(\sigma\nabla u)=0&\text{in }\Omega\\
u=\varphi&\text{on }\partial\Omega
\end{array}\right. 
\end{equation}
and $\tilde{\psi}\in H^1(\Omega)$ is such that
$\tilde{\psi}=\psi$ on $\partial\Omega$ in the trace sense.
We have that $DN(\sigma)$ is a well-defined bounded linear operator,
whose norm is bounded by a constant depending on $N$, $\Omega$, $\lambda_0$, and $\lambda_1$ only, for any $\sigma\in\mathcal{M}(\lambda_0,\lambda_1)$.
Let us notice that, actually, we have $DN(\sigma):H^{1/2}(\partial\Omega)\to H^{-1/2}_{\ast}(\partial\Omega)$. Moreover, since for any constant function $\varphi$ on $\Omega$ we have that $\varphi\in H^{1/2}(\Omega)$ and $DN(\sigma)(\varphi)=0$, no matter what $\sigma$ is, without loss of generality, we actually define
\begin{equation}\label{DNdef}
DN(\sigma):H^{1/2}_{\ast}(\partial\Omega)\to H^{-1/2}_{\ast}(\partial\Omega).\end{equation}

For any conductivity tensor $\sigma$ in $\Omega$, we define its 
Neumann-to-Dirichlet map
$$ND(\sigma):H^{-1/2}_{\ast}(\partial\Omega)\to H^{1/2}_{\ast}(\partial\Omega)$$
where for each $g\in H^{-1/2}_{\ast}(\partial\Omega)$,
$$ND(\sigma)(g)=v|_{\partial\Omega}$$
where $v$ solves
\begin{equation}\label{Neumannproblem1}
\left\{\begin{array}{ll}
\mathrm{div}(\sigma\nabla v)=0&\text{in }\Omega\\
\sigma\nabla v\cdot \nu=g&\text{on }\partial\Omega\\
\int_{\partial \Omega}v=0.&
\end{array}\right. 
\end{equation}
We have that $ND(\sigma)$ is a well-defined bounded linear operator, it is 
the inverse of $DN(\sigma)$ as defined in \eqref{DNdef}, and its
norm is bounded by a constant depending on $N$, $\Omega$, and $\lambda_0$ only, for any $\sigma\in\mathcal{M}(\lambda_0,\lambda_1)$.

We consider the following forward operators 
$$
DN:\mathcal{M}(\lambda_0,\lambda_1)
\to \mathcal{L}(H^{1/2}_{\ast}(\partial \Omega),H^{-1/2}_{\ast}(\partial \Omega))$$
and
$$ND:\mathcal{M}(\lambda_0,\lambda_1)
\to \mathcal{L}(H^{-1/2}_{\ast}(\partial \Omega),H^{1/2}_{\ast}(\partial \Omega)).$$

We can state the inverse conductivity problem in the following way. We wish to determine an unknown conductivity tensor $\sigma$ in $\Omega$ by 
performing electrostatic measurements at the boundary of voltage and current type.
If all boundary measurements are performed, this is equivalent to say that we are measuring either its Dirichlet-to-Neumann map $DN(\sigma)$ or its  Neumann-to-Dirichlet map $ND(\sigma)$. In other words, given either $DN(\sigma)$ or $ND(\sigma)$,
we wish to recover $\sigma$.

Such an inverse problem has a long history, it was in fact proposed by Calder\'on \cite{Cal} in 1980. About uniqueness, there are several result for scalar, that is isotropic, conductivities. In dimension $3$ and higher, already in the 80's, uniqueness was proved in \cite{Koh-Vog84:1,Koh-Vog85} for the determination of the conductivity at the boundary and for the
 analytic case, and then in \cite{Syl-Uhl87} for $C^2$ conductivities.
Slightly later it appeared the first uniqueness result for smooth conductivities in dimension $2$, \cite{Nac}.

Recently, the two dimensional case was completely solved, \cite{Ast-Pai}, for $L^{\infty}$ scalar conductivities. Also for the $N$ dimensional case, with $N\geq 3$, there has been a great improvement. In \cite{Hab-Tat}, the regularity has been reduced to $C^1$ or Lipschitz but close to a constant. The case of general Lipschitz
conductivities is treated in \cite{Car-Rog}. The most general result is the one in \cite{Hab}, where conductivities with unbounded gradient are allowed
and uniqueness is shown for $W^{1,N}$ conductivities, at least for $N=3,4$.

For what concerns anisotropic conductivities, for instance when we consider symmetric conductivity tensors in $\mathcal{M}_{sym}$, uniqueness is never achieved. In fact, let $\varphi:\Omega\to\Omega$ be a bi-Lipschitz mapping, that is a bijective map such that
$\varphi$ and its inverse $\varphi^{-1}$ are Lipschitz functions. Clearly $\varphi$ can be extended to a Lipschitz function defined on $\overline{\Omega}$.
For any $\sigma\in \mathcal{M}_{sym}$ in $\Omega$ and any of these bi-Lipschitz mapping $\varphi$ from $\Omega$ onto itself, we define
the push-forward of the conductivity tensor $\sigma$ by $\varphi$ as
\begin{equation}\label{push-forward}
\varphi_{\ast}(\sigma)(y)=\frac{J(x)\sigma(x)J(x)^T}{|\det J(x)|}\quad\text{for almost any }y\in\Omega
\end{equation}
where $J(x)=J\varphi(x)$ is the Jacobian matrix of $\varphi$ in $x$ and $x=\varphi^{-1}(y)$.
We have that $\varphi_{\ast}(\sigma)\in \mathcal{M}_{sym}$ and that
\begin{equation}\label{nonunique}
DN(\sigma)=DN(\varphi_{\ast}(\sigma))\text{ and }ND(\sigma)=ND(\varphi_{\ast}(\sigma))\quad\text{if }\varphi|_{\partial\Omega}=Id.
\end{equation}
In dimension $N=2$ and for $\Omega$ simply connected, \eqref{push-forward} and
\eqref{nonunique} still hold even if we consider $\varphi:\Omega\to\Omega$ to be a 
quasiconformal mapping. We recall that, for $\Omega\subset \mathbb{R}^2$, simply connected bounded open set with Lipschitz boundary, we say that $\varphi:\Omega\to\Omega$ is a quasiconformal mapping if $\varphi$ is bijective, $\varphi\in W^{1,2}(\Omega)$, and, for some $K\geq 1$, we have
$$\|J\varphi(x)\|^2\leq K\det(J\varphi(x))\quad \text{for a.e. }x\in\Omega.$$

By \eqref{nonunique}, it is immediate to notice that our inverse problem can not have a unique solution if we consider symmetric conductivity tensors. On the other hand, 
in dimension $2$, this is the only obstruction to uniqueness for symmetric conductivity tensors, as proved in \cite{Syl} in the smooth case and in \cite{Ast-Pai-Las} in the general $L^{\infty}$ case.

We summarise these results in the following theorem.

\begin{teo}\label{uniqthm}
Let $\Omega\subset \mathbb{R}^N$, $N=2,3,4$, be a bounded, connected domain with Lipschitz boundary.
Let $\sigma_1$ and $\sigma_2$ belong to $\mathcal{M}_{scal}$.

If $N=3,4$ and $\sigma_1$, $\sigma_2\in W^{1,N}(\Omega)$, then we have, see \textnormal{\cite{Hab}},
$$DN(\sigma_1)=DN(\sigma_2)\text{ or }ND(\sigma_1)=ND(\sigma_2)\text{ implies }\sigma_1=\sigma_2.$$

If $N=2$ and $\Omega$ is simply connected, then we have, see \textnormal{\cite{Ast-Pai}},
$$DN(\sigma_1)=DN(\sigma_2)\text{ or }ND(\sigma_1)=ND(\sigma_2)\text{ implies }\sigma_1=\sigma_2.$$

If $N=2$ and $\Omega$ is simply connected,
for any $\sigma\in \mathcal{M}_{sym}$ we define
\begin{multline*}
\Sigma(\sigma)=\{\sigma_1\in\mathcal{M}_{sym} :\ \sigma_1=\varphi_{\ast}(\sigma)\\\text{ where }\varphi:\Omega\to\Omega\text{ is a quasiconformal mapping and }\varphi|_{\partial\Omega}=Id\}.
\end{multline*}
Then $DN(\sigma)$, or equivalently
$ND(\sigma)$, uniquely determines the class $\Sigma(\sigma)$, see \textnormal{\cite{Ast-Pai-Las}}.
 \end{teo}

\subsection{Variational formulation and ill-posedness}\label{illsec}

In practice, the inverse problem consists in the following. Let $\sigma_0$ be a conductivity tensor in $\Omega$ that we wish to determine. Considering for example the Dirichlet-to-Neumann case, we measure $DN(\sigma_0)$. Since our measurements are obviously noisy, the information that is actually available is a perturbation of $DN(\sigma_0)$, that we may call $\hat{\Lambda}$. Therefore our inverse problem consists in finding a conductivity $\sigma$ such that
$DN(\sigma)=\hat{\Lambda}$. Due to the noise in the measurements this problem may not have any solution. We should therefore solve the problem in a least-square-type way, namely
solve
$$
\min_{\sigma}\|DN(\sigma)-\hat{\Lambda}\|.
$$
The fact that such a minimum problem admits a solution depends on several aspects. In particular it depends on the class of conductivity tensors on which we consider the minimisation and, in part, also on the kind of norm we use to measure the distance between $DN(\sigma)$ and $\hat{\Lambda}$. Next we discuss in details these issues.

Occasionally, we shall use the so-called $H$-convergence. For a definition and its basic properties we refer to \cite{All,Mur-Tar1,Mur-Tar2}. We recall that $G$- or $H$-convergence was shown to be quite useful for the inverse conductivity problem, see for instance \cite{Ale-Cab,Far-Kur-Rui,Ron15}.
Here we just remark a few of its properties. This is a very weak kind of convergence, in fact it is weaker than $L^1_{loc}$ convergence.
For 
symmetric conductivity tensors $H$-convergence reduces to the more usual $G$-convergence.  The most important fact is that
$\mathcal{M}(\lambda_0,\lambda_1)$ is compact with respect to $H$-convergence
and $\mathcal{M}_{sym}(\lambda_0,\lambda_1)$ is also compact with respect to $H$-convergence, or equivalently $G$-convergence. Furthermore, $\mathcal{M}_{scal}(\lambda_0,\lambda_1)$ is not closed with respect to $G$-convergence, actually
any symmetric conductivity tensor is the limit, in the $G$-convergence sense, of scalar conductivities assuming only two different positive values.

We use the following notation. Let $B_1$ and $B_2$ be two Banach spaces such that $B_1\subset H^{1/2}_{\ast}(\partial\Omega)$
and $H^{-1/2}_{\ast}(\partial\Omega)\subset B_2$, with continuous immersions.
Moreover, let $\tilde{B}_1$ and $\tilde{B}_2$ be two Banach spaces such that $\tilde{B}_1\subset H^{-1/2}_{\ast}(\partial\Omega)$
and $H^{1/2}_{\ast}(\partial\Omega)\subset \tilde{B}_2$, with continuous immersions.

We denote with $X$ the space $\mathcal{M}(\lambda_0,\lambda_1)$, $\mathcal{M}_{sym}(\lambda_0,\lambda_1)$, or
$\mathcal{M}_{scal}(\lambda_0,\lambda_1)$. The natural metric on $X$ will be the one induced by the $L^1$ metric.

In the Dirichlet-to-Neumann case, we call $Y=\mathcal{L}(B_1,B_2)$, with the distance induced by its norm, and denote $\Lambda=DN:X\to Y$. 

We speak of the natural norm of the Dirichlet-to-Neumann map when
$B_1=H^{1/2}_{\ast}(\partial\Omega)$ and $B_2=H^{-1/2}_{\ast}(\partial\Omega)$ and we denote it with $\|\cdot\|_{nat}$ or $\|\cdot\|_{H^{1/2},H^{-1/2}}$.
We have a canonical continuous linear map from $\mathcal{L}(H^{1/2}_{\ast}(\partial\Omega),H^{-1/2}_{\ast}(\partial\Omega))$ into $Y$.
If we assume that $B_1$ is dense in $H^{1/2}_{\ast}(\partial\Omega)$, then this map is injective, thus $\mathcal{L}(H^{1/2}_{\ast}(\partial\Omega),H^{-1/2}_{\ast}(\partial\Omega))\subset Y$, with continuous immersion, and, if $y\in Y$ is such that
$\|y\|_Y=0$, then
$y\in \mathcal{L}(H^{1/2}_{\ast}(\partial\Omega),H^{-1/2}_{\ast}(\partial\Omega))$ 
and also $\|y\|_{nat}=0$.

In the Neumann-to-Dirichlet case, we call $Y=\mathcal{L}(\tilde{B}_1,\tilde{B}_2)$, with the distance induced by its norm, and denote $\Lambda=ND:X\to Y$. 

We speak of the natural norm of the Neumann-to-Dirichlet map when
$\tilde{B}_1=H^{-1/2}_{\ast}(\partial\Omega)$ and $\tilde{B}_2=H^{1/2}_{\ast}(\partial\Omega)$
and we denote it with $\|\cdot\|_{nat}$ or $\|\cdot\|_{H^{-1/2},H^{1/2}}$.
We have a canonical continuous linear map from $\mathcal{L}(H^{-1/2}_{\ast}(\partial\Omega),H^{1/2}_{\ast}(\partial\Omega))$ into $Y$.
If we assume that $\tilde{B}_1$ is dense in $H^{-1/2}_{\ast}(\partial\Omega)$, then this map is injective, thus $\mathcal{L}(H^{-1/2}_{\ast}(\partial\Omega),H^{1/2}_{\ast}(\partial\Omega))\subset Y$, with continuous immersion, and, if
if $y\in Y$ is such that
$\|y\|_Y=0$, then
$y\in \mathcal{L}(H^{-1/2}_{\ast}(\partial\Omega),H^{1/2}_{\ast}(\partial\Omega))$ 
and also $\|y\|_{nat}=0$.
Another interesting and useful choice for $\tilde{B}_1$ and $\tilde{B}_2$ is given by
$\tilde{B}_1=\tilde{B}_2=L^2_{\ast}(\partial\Omega)$, see the discussion in \cite{Ron15}, 
and we denote its norm with $\|\cdot\|_{L^2,L^2}$. We remark that
$L^2_{\ast}(\partial\Omega)$ is clearly dense in $H^{-1/2}_{\ast}(\partial\Omega)$.

Let us notice in the following remark that, when we consider the natural norms, then all results related to the Dirichlet-to-Neumann maps may be proved also for the Neumann-to-Dirichlet maps, and viceversa.

\begin{oss}\label{DNvsND}
Let $\sigma_1$, $\sigma_2\in \mathcal{M}(\lambda_0,\lambda_1)$.
Then there exist positive constants $C_1$ and $C_2$, depending on $N$, $\Omega$,
$\lambda_0$, and $\lambda_1$ only, such that
\begin{multline*}
C_1\|ND(\sigma_1)-ND(\sigma_2)\|_{H^{-1/2},H^{1/2}}\leq\|DN(\sigma_1)-DN(\sigma_2)\|_{H^{1/2},H^{-1/2}}\\\leq C_2\|ND(\sigma_1)-ND(\sigma_2)\|_{H^{-1/2},H^{1/2}}.
\end{multline*}

In fact, we have
$$DN(\sigma_1)-DN(\sigma_2)=DN(\sigma_1)(ND(\sigma_2)-ND(\sigma_1))DN(\sigma_2)$$
and the same formula holds if we swap $DN$ with $ND$.
\end{oss}

If we call $\hat{\Lambda}\in Y$ either the measured Dirichlet-to-Neumann map or the measured Neumann-to-Dirichlet map,
then the inverse problem consists in finding $\sigma\in X$ such that
$\Lambda(\sigma)=\hat{\Lambda}$. However, since
$\hat{\Lambda}$ is a measured, therefore noisy, quantity, this problem may not have any solution and we thus solve the problem in a least-square-type way, namely
solve
\begin{equation}\label{minpbmcond}
\min\{\|\Lambda(\sigma)-\hat{\Lambda}\|_{Y}:\ \sigma\in X\}.
\end{equation}

Such a problem always admits a solution either if $X=\mathcal{M}(\lambda_0,\lambda_1)$ or if $X=\mathcal{M}_{sym}(\lambda_0,\lambda_1)$.
In fact the following is proved in \cite{Ron15}.

\begin{prop}\label{Hconv}
Under the previous notation and assumptions, 
let us consider a sequence of conductivity tensors $\{\sigma_n\}_{n\in\mathbb{N}}\subset \mathcal{M}(\lambda_0,\lambda_1)$ and a conductivity tensor $\sigma$ in the same set.

If, as $n\to\infty$, $\sigma_n$ converges to $\sigma$ strongly in $L^1_{loc}$ or in the $H$-convergence sense, then
$$\|\hat{\Lambda}-\Lambda(\sigma)\|_Y\leq\liminf_n \|\hat{\Lambda}-\Lambda(\sigma_n)\|_Y.$$

If $X$ is equal to $\mathcal{M}(\lambda_0,\lambda_1)$ or to $\mathcal{M}_{sym}(\lambda_0,\lambda_1)$,
by compactness of $X$ with respect to $H$-convergence, we deduce that
\eqref{minpbmcond} admits a solution.
\end{prop}

On the other hand, if $X$ is equal to $\mathcal{M}_{scal}(\lambda_0,\lambda_1)$ then
\eqref{minpbmcond} may fail to have a solution as we shall see later on in Example~\ref{noexistenceexam}.

We notice that Proposition~\ref{Hconv} contains a lower semicontinuity result. For certain application, instead, continuity is needed. For our purposes it will be enough the following result, proved in \cite{Ale-Cab}.

\begin{prop}\label{Gcont}
Under the previous notation and assumptions, 
let us consider a sequence of symmetric conductivity tensors $\{\sigma_n\}_{n\in\mathbb{N}}\subset \mathcal{M}_{sym}(\lambda_0,\lambda_1)$ and a conductivity tensor $\sigma$ in the same set. We assume that for some $\Omega'$ compactly contained in $\Omega$ we have $\sigma_n=\sigma$ almost everywhere in $\Omega\backslash\Omega'$ for any $n\in\mathbb{N}$.

If, as $n\to\infty$, $\sigma_n$ converges to $\sigma$ strongly in $L^1_{loc}$ or in the $G$-convergence sense, then
$$\lim_n\|\Lambda(\sigma)-\Lambda(\sigma_n)\|_{nat}=0$$
as well as in $\|\cdot\|_Y$ for any $Y$ as above. 
\end{prop}

We notice that a certain control of the conductivity tensors near the boundary is indeed needed, see \cite[Theorem~4.9]{Far-Kur-Rui}. In the same paper a more general and essentially optimal version of Proposition~\ref{Gcont} is proved, see \cite[Theorem~1.1]{Far-Kur-Rui}.

Proposition~\ref{Gcont} is enough to show that \eqref{minpbmcond} may fail to have a solution if 
$X=\mathcal{M}_{scal}(\lambda_0,\lambda_1)$. We slightly generalise \cite[Example~3.4]{Ron15}, which is based on a nice remark by Giovanni, which is presented in \cite{Ron15} as Example~4.4.
This generalisation shows that existence may fail for both the Dirichlet-to-Neumann and Neumann-to-Dirichlet case and for the natural norms, as well as for any $\|\cdot\|_Y$ with $Y$ as above, if $B_1$ is dense in $H^{1/2}_{\ast}(\partial\Omega)$ or $\tilde{B}_1$ is dense in $H^{-1/2}_{\ast}(\partial\Omega)$, respectively.
It firstly appeared in \cite{Des}, and we present its proof here for the convenience of the reader.

\begin{exam}\label{noexistenceexam}
Let $\Omega=B_1(0)\subset\mathbb{R}^2$.
Under the previous notation and assumptions, let us assume that
$B_1$ is dense in $H^{1/2}_{\ast}(\partial\Omega)$ or $\tilde{B}_1$ is dense in $H^{-1/2}_{\ast}(\partial\Omega)$, respectively.

Let $a>0$ be a positive constant with $a\neq 1$. We define the conductivity
tensor $\tilde{\sigma}\in \mathcal{M}_{sym}(\lambda_0,\lambda_1)\backslash\mathcal{M}_{scal}(\lambda_0,\lambda_1)$ in $B_1(0)\subset \mathbb{R}^2$ as follows
\begin{equation}
\tilde{\sigma}=\left\{\begin{array}{ll}
I_2 &\text{in }B_1(0)\backslash B_{1/2}(0)\\
\left[\begin{matrix}a& 0\\0 &a^{-1}\end{matrix}\right]&\text{in }B_{1/2}(0).\\
\end{array}\right.
\end{equation}

Let us set $\hat{\Lambda}=\Lambda(\tilde{\sigma})$. There exist $0<\lambda_0<\lambda_1$ such that the minimum problem
$$\min_{\sigma\in\mathcal{M}_{scal}(\lambda_0,\lambda_1)}
\|\Lambda(\tilde{\sigma})-\Lambda(\sigma)\|_Y
$$
does not have any solution, for any $Y$ as above, thus including the natural norms.
\end{exam}
 
\proof{.} The crucial point is the following. By density of scalar conductivities
inside symmetric conductivity tensors that follows by the results in \cite{Mur-Tar2}, see
\cite[Proposition~2.2]{Ron15} for a convenient version, we can find $0<\lambda_0<\lambda_1$ and $\{\sigma_n\}_{n\in\mathbb{N}}\subset \mathcal{M}_{scal}(\lambda_0,\lambda_1)$ such that
$\sigma_n$ $G$-converges to $\tilde{\sigma}$ as $n\to\infty$ and
$\sigma_n=I_2$ in $B_1(0)\backslash B_{1/2}(0)$ for any $n\in\mathbb{N}$.
Therefore, by Proposition~\ref{Gcont}, we immediately conclude that 
$$\inf_{\sigma\in\mathcal{M}_{scal}(\lambda_0,\lambda_1)}
\|\Lambda(\tilde{\sigma})-\Lambda(\sigma)\|_Y=0.$$

In order for a minimiser to exist, then we need to find a scalar conductivity
$\hat{\sigma}$ such that
$\|\Lambda(\hat{\sigma})-\Lambda(\tilde{\sigma})\|_Y=0$, hence, by our density assumptions, such that
 $\Lambda(\hat{\sigma})=\Lambda(\tilde{\sigma})$.
By the main result of \cite{Ast-Pai-Las}, recalled in Theorem~\ref{uniqthm},
there exists a quasiconformal mapping $\varphi:B_1(0)\to B_1(0)$ such that $\varphi|_{\partial B_1(0)}=Id$ and
$\varphi_{\ast}(\hat{\sigma})=\tilde{\sigma}$. We recall that actually
$\varphi:\overline{B_1(0)}\to\overline{B_1(0)}$, it is continuous, bijective and its inverse is continuous as well.
We assume that $\hat{\sigma}(x)=s(x)I_2$, $x\in
B_1(0)$, with $s\in L^{\infty}(B_1(0))$ and bounded away from $0$. Then $\varphi_{\ast}(\hat{\sigma})=\tilde{\sigma}$ means that for almost any $y\in B_1(0)$ we have
\begin{multline*}
\tilde{\sigma}(y)=\frac{J(x)(s(x) I_2)J(x)^{T}}{|\det J(x)|}\\
=\frac{s(x)}{|\det J(x)|}\left[\begin{matrix}|\nabla \varphi_1(x)|^2& \nabla\varphi_1(x)\cdot\nabla\varphi_2(x)\\ \nabla\varphi_1(x)\cdot\nabla\varphi_2(x) &|\nabla \varphi_2(x)|^2\end{matrix}\right],
\end{multline*}
where $\varphi=(\varphi_1,\varphi_2)$,  $J(x)$ is the Jacobian matrix of
$\varphi$ in $x$, and $x=\varphi^{-1}(y)$.
Since $\det(\tilde{\sigma}(y))=1$ for almost any $y\in B_1(0)$,
we conclude that, for almost any $x\in B_1(0)$, $s(x)=1$, that is
$\hat{\sigma}\equiv I_2$ in $B_1(0)$. We also note that, since $\varphi$ is quasiconformal, then $\det J(x)>0$ for almost any $x\in B_1(0)$.

By the structure of $\tilde{\sigma}$, we infer that
for almost any $x\in B_1(0)$ we have
$\nabla \varphi_2(x)= \lambda(x) \left[\begin{smallmatrix}0 & -1\\1 & 0
\end{smallmatrix}\right]\nabla \varphi_1(x)$ with $\lambda(x)>0$, since $\det J(x)>0$, satisfying the following
$$\lambda=\left\{\begin{array}{ll}
1&\text{in }D=\varphi^{-1}(B_1(0)\backslash B_{1/2}(0))\\
a^{-1}&\text{in }D_1=B_1(0)\backslash D=\varphi^{-1}(B_{1/2}(0)).
\end{array}\right.
$$
We conclude that
$$\Delta \varphi_1=\Delta\varphi_2=0\quad\text{in }D\text{ and in }D_1.$$
More precisely, we have that $\varphi_1+\rmi\varphi_2$ is holomorphic in $D$.
Since $\varphi_1(x_1,x_2)=x_1$ and $\varphi_2(x_1,x_2)=x_2$ on $\partial B_1(0)$, by the unique continuation from Cauchy data, we infer that $\varphi=Id$ in $D$ as well.
Therefore $B_1(0)\backslash B_{1/2}(0)=\varphi(D)=D$. We conclude that
$D_1=B_{1/2}(0)$, $\varphi_1$ and $\varphi_2$ are harmonic in $B_{1/2}(0)$,
and  $\varphi_1(x_1,x_2)=x_1$ and $\varphi_2(x_1,x_2)=x_2$ on $\partial B_{1/2}(0)$.
We immediately conclude that $\varphi=Id$ on the whole $B_1(0)$ and we obtain a contradiction.\cvd

\bigskip

If we have no control on the conductivity tensors near the boundary, then continuity of our forward operators may be achieved by suitably choosing the spaces $B_1$, $B_2$, and $\tilde{B}_1$, $\tilde{B}_2$, that is by changing
the distance, thus the space $Y$, with respect to which we measure the error on our measurements.
Namely we have the following results, see \cite{Ron15}.

\begin{prop}\label{Hconvcont}
Under the previous notation and assumptions,
there exists $Q_1>2$, depending on $N$, $\Omega$, $\lambda_0$, and $\lambda_1$ only,
such that
the following holds for any 
$2<p<Q_1$.

In the Dirichlet-to-Neumann case, we assume that
$B_1\subset W^{1-1/p,p}_{\ast}(\partial \Omega)$,
with continuous immersion.

In the Neumann-to-Dirichlet case, we assume that
$\tilde{B}_1\subset (W^{1-1/p',p'}(\partial \Omega))'_{\ast}$, with continuous immersion, where $(W^{1-1/p',p'}(\partial \Omega))'_{\ast}$ is the subspace of $g$ belonging to the
dual of $W^{1-1/p',p'}(\partial \Omega)$ such that $\langle g,1\rangle =0$.

Then $\Lambda$ is H\"older continuous with respect to the $L^1$ distance in $\mathcal{M}(\lambda_0,\lambda_1)$ and the distance $d$ on $Y$ given by its norm.
The H\"older exponent $\beta$ is equal to $(p-2)/(2p)$.
\end{prop}

A particularly interesting case for Neumann-to-Dirichlet maps is to choose $\tilde{B}_1=\tilde{B}_2=L^2_{\ast}(\partial\Omega)$ since 
$L^2(\partial\Omega)$ is contained in the dual of $W^{1-1/p',p'}(\partial \Omega)$ for some $p$, $2<p<Q_1$, with $p$ close enough to $2$,
and $H^{1/2}_{\ast}(\partial\Omega)\subset L^2_{\ast}(\partial\Omega)$, with continuous immersions. Moreover, $L^2_{\ast}(\partial\Omega)$ is dense in $H^{-1/2}_{\ast}(\partial\Omega)$. In this case we also have continuity with respect to $H$-convergence, see again \cite{Ron15}.

\begin{prop}\label{Hconvcont2}
Under the previous notation and assumptions,
let us consider a sequence of conductivity tensors $\{\sigma_n\}_{n\in\mathbb{N}}\subset \mathcal{M}(\lambda_0,\lambda_1)$ and a conductivity tensor $\sigma$ in the same set.

If, as $n\to\infty$, $\sigma_n$ converges to $\sigma$ strongly in $L^1_{loc}$ or in the $H$-convergence sense, then
$$\lim_n\|ND(\sigma)-ND(\sigma_n)\|_{L^2,L^2}=0.$$
\end{prop}

Let us consider that $\sigma_0\in X$ is the conductivity tensor in $\Omega$ that we wish to determine. Given the noise level $\varepsilon>0$, our measurement is given by $\hat{\Lambda}_{\varepsilon}\in Y$, satisfying
\begin{equation}\label{noiselevel}
\|\hat{\Lambda}_{\varepsilon}-\Lambda(\sigma_0)\|_Y\leq\varepsilon.
\end{equation}
For consistency, we call $\hat{\Lambda}_0=\Lambda(\sigma_0)$.
Assume that our minimisation problem
\begin{equation}\label{minpbmcond2}
\min\{\|\Lambda(\sigma)-\hat{\Lambda}_{\varepsilon}\|_{Y}:\ \sigma\in X\}
\end{equation}
admits a solution and let us call $\tilde{\sigma}_{\varepsilon}$ a minimiser for
\eqref{minpbmcond2}. The main question is whether 
$\tilde{\sigma}_{\varepsilon}$ is a good approximation of the looked for conductivity tensor $\sigma_0$, namely we ask whether $\lim_{\varepsilon\to 0^+}\tilde{\sigma}_{\varepsilon}=\sigma_0$, where the limit is to be intended in a suitable sense.
Unfortunately this may not be true, in fact our inverse problem is ill-posed, that is, we have no stability. There are two serious obstructions to stability. In the anisotropic case, that is, when $X=\mathcal{M}_{sym}(\lambda_0,\lambda_1)$, for instance, the obstruction is due to invariance by changes of coordinates that keep fixed the boundary.
In the isotropic case, that is, when $X=\mathcal{M}_{scal}(\lambda_0,\lambda_1)$, the obstruction is due to the fact that this class is not closed with respect to $G$-convergence.

Let us illustrate these difficulties in the following three examples.

\begin{exam}\label{examanisot}
Let $\Omega=B_1(0)\subset\mathbb{R}^2$.
Let $X=\mathcal{M}_{sym}(\lambda_0,\lambda_1)$, for some $0<\lambda_0<1<\lambda_1$ to be fixed later.
We set $\tilde{\sigma}_0\equiv I_2$ in $B_1(0)\subset\mathbb{R}^2$.
We fix a $C^1$ diffeomorphism $\varphi:B_1(0)\to B_1(0)$ such that 
$\varphi$ is identically equal to the identity in $B_1(0)\backslash B_{1/2}(0)$. We call
$\sigma_0=\varphi_{\ast}(\tilde{\sigma}_0)$ and we assume that $\sigma_0$ is the conductivity tensor to be recovered. We notice that, if $\varphi$ is not trivial, we have that
$\sigma_0\neq\tilde{\sigma}_0$.

Let $\tilde{\sigma}_{\varepsilon}$, $0<\varepsilon\leq \varepsilon_0$, be a scalar conductivity satisfying
$\|\tilde{\sigma}_{\varepsilon}-\tilde{\sigma}_0\|_{L^{\infty}(\Omega)}\leq \varepsilon$.
We notice that, choosing in a suitable way $\lambda_0$ and $\lambda_1$,  
we have $\sigma_0\in X$, and,
for any 
$0\leq\varepsilon\leq \varepsilon_0$, also 
$\tilde{\sigma}_{\varepsilon}\in X$.

We notice that $\Lambda(\sigma_0)=\Lambda(\tilde{\sigma}_0)$ and, for some constant $C$, depending on $\lambda_0$, $\lambda_1$, and $Y$ only, we have
$$
\|\Lambda(\tilde{\sigma}_{\varepsilon})-\Lambda(\sigma_0)\|_Y\leq C\varepsilon.
$$

If, for any 
$0<\varepsilon\leq \varepsilon_0$, we assume that
$\hat{\Lambda}_{\varepsilon}=\Lambda(\tilde{\sigma}_{\varepsilon})$,
then,
unfortunately,  we have
$$\tilde{\sigma}_{\varepsilon}\in\argmin_{\sigma\in X}\|\Lambda(\sigma)-\hat{\Lambda}_{\varepsilon}\|_{Y},
$$
and, obviously, for any sequence $\{\varepsilon_n\}_{n\in\mathbb{N}}\subset (0,\varepsilon_0]$ such that $\lim_n\varepsilon_n=0$,
$\tilde{\sigma}_{\varepsilon_n}$
does not converge, not even in the $G$-convergence sense or in the weak $L^1$ sense, to $\sigma_0$.
\end{exam}

In dimension $2$, in \cite{Far-Kur-Rui}, it has been proved that this is the only obstruction in the symmetric conductivity tensor case, if we consider the natural norms. Namely, from \cite[Theorem~1.3]{Far-Kur-Rui}, we can immediately deduce the following.

\begin{prop}\label{weakconprop}
Let $N=2$ and let $\Omega\subset \mathbb{R}^2$ be a bounded, simply connected open set with Lipschitz boundary.
Let $\sigma_0\in X=\mathcal{M}_{sym}(\lambda_0,\lambda_1)$.
We pick either $Y=\mathcal{L}(H^{1/2}_{\ast}(\partial\Omega),H^{-1/2}_{\ast}(\partial\Omega))$, for the Dirichlet-to-Neumman case, or $Y=\mathcal{L}(H^{-1/2}_{\ast}(\partial\Omega),H^{1/2}_{\ast}(\partial\Omega))$, for the Neumann-to-Dirichlet case, respectively.
For any $n\in\mathbb{N}$, let $\hat{\Lambda}_n\in Y$ be such that
$$\|\hat{\Lambda}_n-\Lambda(\sigma_0)\|_{Y}=\|\hat{\Lambda}_n-\Lambda(\sigma_0)\|_{nat}\to 0\quad\text{as }n\to\infty.$$
Let $\sigma_n\in\argmin_{\sigma\in X} \|\hat{\Lambda}_n-\Lambda(\sigma_0)\|_Y$, $n\in\mathbb{N}$. Then, for any $n\in\mathbb{N}$, there exists a quasiconfomal mapping $\varphi_n:\Omega\to\Omega$ such that $(\varphi_n)|_{\partial\Omega}=Id$ and
$$(\varphi_n)_{\ast}(\sigma_n)\to\sigma_0\quad\text{as }n\to\infty$$
in the $G$-convergence sense.
\end{prop}

\proof{.} We have that
\begin{multline*}
\|\Lambda(\sigma_n)-\Lambda(\sigma_0)\|_Y\leq \|\hat{\Lambda}_n-\Lambda(\sigma_n)\|_Y+ \|\Lambda(\sigma_0)-\hat{\Lambda}_n\|_Y\\
\leq 2\|\Lambda(\sigma_0)-\hat{\Lambda}_n\|_Y\to 0\quad\text{as }n\to\infty.
\end{multline*}
Then the conclusion follows by \cite[Theorem~1.3]{Far-Kur-Rui}.\cvd

\bigskip

We notice that the kind of convergence we have in Proposition~\ref{weakconprop} is really weak, in several respects. First, it is only up to a change of variables, second it is in the sense of $G$-convergence, only. We recall that $G$-convergence does not imply
convergence not even in the weak $L^1$ sense. In fact, let us consider the following example. Let $D$ be an open set such that $D\subset Q=(0,1)^N\subset\mathbb{R}^N$, $N\geq 2$, and let us consider, for two given constants $0<a<b$,
\begin{equation}\label{periodic}
\sigma=\left\{\begin{array}{ll}
a&\text{in }D\\
b&\text{in }Q\backslash D.
\end{array}\right.
\end{equation}
We also assume that $D$ and $Q\backslash D$ have positive measure.
Then we have
$$a<m_h=\left(\int_Q \sigma^{-1}\right)^{-1}<\int_Q\sigma=m<b $$
where $m_h$ is the so-called harmonic mean of $\sigma$ on $Q$ and $m$ is the usual mean of $\sigma$ on $Q$.

We extend $\sigma$ all over $\mathbb{R}^N$ by periodicity and define, for any $\varepsilon>0$,
$$\sigma_{\varepsilon}(x)=\sigma(x/\varepsilon)I_N,\quad x\in\mathbb{R}^N.$$
Given $\Omega$ a bounded connected open set with Lipschitz boundary,
it is a classical fact in homogenisation theory that in $\Omega$
$$\sigma_{\varepsilon}\text{ $G$-converges to }\sigma_{hom}\quad\text{as }\varepsilon\to 0^+$$ 
where $\sigma_{hom}$ is a constant symmetric matrix satisfying
$$m_hI_N\leq \sigma_{hom}<mI_N.$$
On the other hand, $\sigma_{\varepsilon}$ converges to $a_mI_N$ in the weak$^{\ast}$ $L^{\infty}(\Omega)$ sense, therefore also weakly in $L^1(\Omega)$.

Moreover, if $N=2$ and
 \begin{equation}\label{Ddefin}
D=\{(x_1,x_2)\in Q:\ (x_1-1/2)(x_2-1/2)>0\},
 \end{equation}
then $\sigma_{hom}$ can be computed explicitly and we have
that $\sigma_{hom}=\sqrt{ab}I_2$. 

Instead, if $N=2$ and
\begin{equation}\label{Ddefinbis}
D=\{(x_1,x_2)\in Q:\ (x_1-1/2)>0\},
 \end{equation}
then also in this case $\sigma_{hom}$ can be computed explicitly and we have
that
$$\sigma_{hom}=\left[\begin{matrix}m_h &0 \\ 0 &m \end{matrix}\right].$$

These explicit formulas are the bases for the next examples. The next one
was introduced in \cite{Ale-Cab} and we state it here. This and the next example show that in the scalar case, when, at least in dimension $2$, uniqueness is not an issue, instability phenomena may occur, no matter what we choose as $Y$.

\begin{exam}\label{AleCabexam}
Let $\Omega=B_1(0)\subset\mathbb{R}^2$.
Let $X=\mathcal{M}_{scal}(\lambda_0,\lambda_1)$, for some $0<\lambda_0<1<\lambda_1$ to be fixed later.
Let us assume that
$B_1$ is dense in $H^{1/2}_{\ast}(\partial\Omega)$ or $\tilde{B}_1$ is dense in $H^{-1/2}_{\ast}(\partial\Omega)$, respectively.

We fix $N=2$ and two positive constants $0<a<b$. We take $Q=(0,1)^2$ and $D$ as in \eqref{Ddefin}.
We call
$$\sigma_0=\left\{\begin{array}{ll}I_2 &\text{in }B_1(0)\backslash B_{1/2}(0)\\
\sqrt{ab} &\text{in }B_{1/2}(0).
\end{array}
\right.
$$
We define $\sigma$ as in \eqref{periodic}, we extend it by periodicity all over $\mathbb{R}^2$, and define, for any $\varepsilon$, $0<\varepsilon\leq 1/2$,
$$\sigma_{\varepsilon}=\left\{\begin{array}{ll}I_2 &\text{in }B_1(0)\backslash B_{1/2}(0)\\
 \sigma(x/\varepsilon)I_N&\text{if }x\in B_{1/2}(0).
\end{array}
\right.
$$

We have that $\sigma_{\varepsilon}$ $G$-converges to $\sigma_0$
as $\varepsilon\to 0^+$, therefore, by Proposition~\ref{Gcont}, we immediately conclude that 
$$\|\Lambda(\tilde{\sigma}_{\varepsilon})-\Lambda(\sigma_0)\|_Y\to
0\quad\text{as }\varepsilon\to 0^+.$$

Therefore, if $\sigma_0$ is the conductivity to be determined, and our measured data
are $\hat{\Lambda}_{\varepsilon}=\Lambda(\sigma_{\varepsilon})$, for any $\varepsilon\in(0,1/2]$, then we have that
$$\{\sigma_{\varepsilon}\}=\argmin_{\sigma\in \mathcal{M}_{scal}(\lambda_0,\lambda_1) }\|\Lambda(\sigma)-\hat{\Lambda}_{\varepsilon}\|_{Y}.$$

On the other hand, we have that $\sigma_{\varepsilon}$ converges to $mI_N$
in 
the weak$^{\ast}$ $L^{\infty}(\Omega)$ sense, therefore also weakly in $L^1(\Omega)$.
Since $\sqrt{ab}<m$, we obtain that, as $\varepsilon\to 0^+$, $\sigma_{\varepsilon}$ does not converge to $\sigma_0$ even weakly in $L^1(\Omega)$, but only $G$-converges to $\sigma_0$.
\end{exam}

The third and final example, inspired by the one in \cite{Ale-Cab} we just presented,
shows that even $G$-convergence may not be guaranteed. 

\begin{exam}\label{examanisot2}
Let $\Omega=B_1(0)\subset\mathbb{R}^2$.
Let $X=\mathcal{M}_{scal}(\lambda_0,\lambda_1)$, for some $0<\lambda_0<1<\lambda_1$ to be fixed later.
Let us assume that
$B_1$ is dense in $H^{1/2}_{\ast}(\partial\Omega)$ or $\tilde{B}_1$ is dense in $H^{-1/2}_{\ast}(\partial\Omega)$, respectively.

We set the conductivity to be determined as $\sigma_0\equiv I_2$ in $B_1(0)\subset\mathbb{R}^2$.
Let $\varphi:B_1(0)\to B_1(0)$ be a $C^1$ diffeomorphism such that 
$\varphi\equiv Id$ in $B_1(0)\backslash B_{1/2}(0)$ and
$\varphi(x_1,x_2)=(x_1/2,x_2)$ on $B_{1/4}(0)$.
We call
$\tilde{\sigma}_0=\varphi_{\ast}(\sigma_0)$. We have that
$\tilde{\sigma}_0\neq \sigma_0$. In particular, 
$\tilde{\sigma}_0=I_2$ in $B_1(0)\backslash B_{1/2}(0)$ and
$\tilde{\sigma}_0(y)=\left[\begin{smallmatrix}1/2 & 0\\ 0 & 2
\end{smallmatrix}\right]$ for any $y\in B_{1/8}(0)$.

We pick $Q=(0,1)^2$, $D$ as in \eqref{Ddefinbis}, and $\sigma$ as in \eqref{periodic}, with
$a=2-\sqrt{3}$ and $b=2+\sqrt{3}$, so that $\sigma_{hom}=\left[\begin{smallmatrix}1/2 &0 \\ 0 &2 \end{smallmatrix}\right]$.

Then, again by density of scalar conductivities
inside symmetric conductivity tensors,
we can find $0<\lambda_0<\lambda_1$ and $\{\tilde{\sigma}_n\}_{n\in\mathbb{N}}\subset \mathcal{M}_{scal}(\lambda_0,\lambda_1)$ such that
$\tilde{\sigma}_n$ $G$-converges to $\tilde{\sigma}_0$ as $n\to\infty$ and such that,
for any $n\in\mathbb{N}$,
$\tilde{\sigma}_n=I_2$ in $B_1(0)\backslash B_{1/2}(0)$ and
$$\tilde{\sigma}_n(y)=\sigma(ny)\quad\text{for any }y\in B_{1/8}(0),$$
where as usual $\sigma$ is extended by periodicity all over $\mathbb{R}^2$.

We notice that, as $n\to\infty$, in $B_{1/8}(0)$, $\tilde{\sigma}_n$ converges to $2I_N$ in the weak$^{\ast}$ $L^{\infty}$ sense, hence also weakly in $L^1(B_{1/8}(0))$. Therefore, $\sigma_n$ can not converge, not even up to subsequences, to $\sigma_0$, not even weakly in $L^1(B_1(0))$.

By Proposition~\ref{Gcont}, we immediately conclude that 
$$\|\Lambda(\tilde{\sigma}_n)-\Lambda(\sigma_0)\|_Y=
\|\Lambda(\tilde{\sigma}_n)-\Lambda(\tilde{\sigma}_0)\|_Y\to
0\quad\text{as }n\to\infty.$$
If we pick as our measured data $\hat{\Lambda}_n=\Lambda(\tilde{\sigma}_n)$, for any $n\in\mathbb{N}$, then we have that
$$\{\tilde{\sigma}_n\}=\argmin_{\sigma\in \mathcal{M}_{scal}(\lambda_0,\lambda_1) }\|\Lambda(\sigma)-\hat{\Lambda}_n\|_{Y}.$$
Then we have that, as $n\to\infty$, $\tilde{\sigma}_n$ can not converge, not even up to subsequences, to the looked for scalar conductivity $\sigma_0$ either in the $G$-convergence sense or locally weakly in $L^1$, hence, a fortiori, in the $L^1_{loc}$ sense as well.
\end{exam}

\subsection{Regularisation}\label{regularsubs}

The issues for this inverse problem previous highlighted, in particular the ill-posedness, lead naturally to consider a suitable regularisation of the minimisation problem \eqref{minpbmcond}. To fix the ideas we consider a regularisation \`a la  Tikhonov. For a general introduction to Tikhonov regularisation, we refer for instance to \cite{Eng-et-al}. Here we are interested in the case of nonsmooth and possibly discontinuous unknown conductivity tensors, therefore we shall follow the approach developed in \cite{Ron08}. We notice that, in the smooth case, the general theory for convergence of Tikhonov regularised solutions for nonlinear operators, as it was developed in \cite{Eng-Kun-Neu}, see also \cite{Eng-et-al}, may be used and leads also to convergence estimates. For example, for the electrical impedance tomography this approach was used in \cite{LMP}, see also \cite{Ji-Ma}.

Instead, in the nonsmooth case,
our starting point is the regularisation strategy proved in \cite{Ron08}, which we recall now. The key ingredient is $\Gamma$-convergence, see \cite{DM} for a detailed introduction.
Here we just recall the definition and basic properties of
$\Gamma$-convergence.

Let $(X,d)$ be a metric space. Then a sequence
$F_n:X\to [-\infty,+\infty]$, $n\in\mathbb{N}$, $\Gamma$-converges as $n\to\infty$
to a function $F:X\to [-\infty,+\infty]$ if for every $x\in X$ we have
\begin{align}\label{liminf}
&\text{for every sequence $\{x_n\}_{n\in\mathbb{N}}$ converging to $x$ we have}\\
&\hspace{2cm}F(x)\leq \liminf_n F_n(x_n);\nonumber\\
\label{limsup}
& \text{there exists a sequence $\{x_n\}_{n\in\mathbb{N}}$ converging to $x$ such that}\\
&\hspace{2cm}F(x)=\lim_n F_n(x_n).\nonumber
\end{align}
The function $F$ will be called the $\Gamma$-limit of the sequence $\{F_n\}_{n\in\mathbb{N}}$ as
$n\to\infty$ with respect to the metric $d$ and we denote it by
$F=\Gamma\textrm{-}\!\lim_n F_n$.
We recall that condition \eqref{liminf} above is usually called the $\Gamma$-liminf inequality, whereas condition \eqref{limsup} is usually referred to as the existence of a recovery sequence.

We say that the functionals $F_n$, $n\in\mathbb{N}$, are equicoercive
if there exists a compact set
$K\subset X$ such that $\inf_K F_n=\inf_X F_n$ for any $n\in\mathbb{N}$.

The following theorem, usually known as the Fundamental Theorem of
$\Gamma$-convergence, illustrates the motivations for the definition
of such a kind of convergence.

\begin{teo}\label{fundthm}
Let $(X,d)$ be a metric space and let $F_n:X\to [-\infty,+\infty]$,
$n\in\mathbb{N}$, 
be a sequence of functions defined on $X$. If the functionals $F_n$, $n\in\mathbb{N}$, are equicoercive
and
$F=\Gamma\textrm{-}\!\lim_n F_n$, then $F$ admits a minimum over $X$ and we
have
$$\min_X F=\lim_n\inf_X F_n.$$
Furthermore, if $\{x_n\}_{n\in\mathbb{N}}$ is a sequence of points in $X$ which
converges to a point $x\in X$ and
satisfies $\lim_n F_n(x_n)=\lim_n\inf_X F_n$, then $x$ is a minimum
point
for $F$.
\end{teo}

The definition of $\Gamma$-convergence may be extended in a natural way
to families depending on a continuous parameter.
The family of functions $F_{\varepsilon}$, defined for every
$\varepsilon>0$, $\Gamma$-converges to a function $F$ as
$\varepsilon\to 0^+$ if for every sequence $\{\varepsilon_n\}_{n\in\mathbb{N}}$ of positive numbers 
converging to $0$ as $n\to\infty$, we have $F=\Gamma\textrm{-}\!\lim_n F_{\varepsilon_n}$.

We begin with an abstract framework. We consider two metric spaces
$(X,d_X)$ and $(Y,d_Y)$ and a continuous function $\Lambda:X\to Y$.
We also fix $x_0\in X$ and $\Lambda_0=\Lambda(x_0)\in Y$.

For any $\varepsilon>0$, we consider a perturbation of $\Lambda_0$ given by
 $\Lambda_{\varepsilon}\in Y$ such that
$d_Y(\Lambda_{\varepsilon},\Lambda_0)\leq \varepsilon$. Here, and in the sequel, $\varepsilon$ plays the role of the noise level.

A function $R:X\to \mathbb{R}\cup\{+\infty\}$ is called a \emph{regularisation operator} for
the metric space $X$ if $R\not\equiv +\infty$ and, with respect to the metric induced by $d_X$,
$R$ is a lower semicontinuous function such that
for any constant $C>0$
the set
$\{x\in X:\ R(x)\leq C\}$ is a compact subset of $X$.

We consider the following regularised minimum problem, for some $\varepsilon>0$,
\begin{equation}\label{minpbm}
\min_{x\in X}(d_Y(\Lambda(x),\Lambda_{\varepsilon}))^{\alpha}+\tilde{a}R(x)
\end{equation}
where $\tilde{a}>0$ is the \emph{regularisation parameter} and $\alpha$ is a positive parameter. In order to make the regularisation meaningful, we need to choose the
regularisation parameter in terms of the noise level $\varepsilon$, namely we choose
$\tilde{a}=\tilde{a}(\varepsilon)$. A solution to \eqref{minpbm} will be called a \emph{regularised solution}.
To fix the ideas,
given $\varepsilon_0>0$, we assume that for any $\varepsilon$, $0<\varepsilon\leq\varepsilon_0$, $\tilde{a}(\varepsilon)=\tilde{a}\varepsilon^\gamma$, for some positive constants $\tilde{a}$ and $\gamma$. By a simple rescaling argument the minimisation problem \eqref{minpbm}
is equivalent to solve
\begin{equation}\label{minpbm0}
\min_{x\in X} F_{\varepsilon}(x)
\end{equation}
where $F_{\varepsilon}:X\to \mathbb{R}\cup\{+\infty\}$ is defined as follows 
\begin{equation}\label{Fdef}
F_{\varepsilon}(x)=
\frac{(d_Y(\Lambda(x),\Lambda_{\varepsilon}))^{\alpha}}{\varepsilon^{\gamma}}+\tilde{a}R(x)
\quad\text{for any }x\in X.
\end{equation}

We also define $F_0:X\to \mathbb{R}\cup\{+\infty\}$ as follows
\begin{equation}\label{defF0}
F_0(x)=\left\{\begin{array}{ll}
\tilde{a}R(x) &\text{if }\Lambda(x)=\Lambda(x_0)=\Lambda_0\text{ in }Y\\
+\infty &\text{otherwise}
\end{array}
\right.
\end{equation}
 for any $x\in X$.

The following result is proved in \cite{Ron08}, by exploiting $\Gamma$-convergence techniques.

\begin{teo}\label{mainteo}
Let $\Lambda$ be continuous and $R$ be a regularisation operator for $X$. Let us also assume that
$R(x_0)<+\infty$ and $\gamma <\alpha$.

Then  we have that there exists $\min_X F_{\varepsilon}$, for any $\varepsilon$, $0\leq\varepsilon\leq \varepsilon_0$, and
$$\min_X F_0=\lim_{\varepsilon\to 0^+}\min_X F_{\varepsilon}<+\infty.$$

Let $\{\tilde{x}_{\varepsilon}\}_{0<\varepsilon\leq\varepsilon_0}$ satisfy $\lim_{\varepsilon\to 0^+} F_{\varepsilon}(\tilde{x}_{\varepsilon})=
\lim_{\varepsilon\to 0^+}\min_X F_{\varepsilon}$ \textnormal{(}for example
we may pick as $\{\tilde{x}_{\varepsilon}\}_{0<\varepsilon\leq\varepsilon_0}$
a family $\{x_{\varepsilon}\}_{0<\varepsilon\leq\varepsilon_0}$ of minimisers of $F_{\varepsilon}$\textnormal{)}.

Let $\{\varepsilon_n\}_{n\in\mathbb{N}}$ be a sequence of positive numbers converging to $0$ as $n\to\infty$.
Then, up to a subsequence, $\tilde{x}_{\varepsilon_n}$ converges to a point $\tilde{x}\in X$
such that $\tilde{x}$ is a minimiser of $F_0$, that is, in particular, $\Lambda(\tilde{x})=\Lambda(x_0)$ in $Y$
and $R(\tilde{x})=\min\{R(x):\ x\in X\text{ such that }\Lambda(x)=\Lambda(x_0)\text{ in }Y\}$.

Furthermore, if $F_0$ admits a unique minimiser $\tilde{x}$, then
we have that
\begin{equation}\label{conv1}
\lim_{\varepsilon\to 0^+}\tilde{x}_{\varepsilon}=\tilde{x}.
\end{equation}

Finally, if on the set $\{x\in X:\ R(x)<+\infty\}$ the map $\Lambda$ is injective, then
we have
$$\lim_{\varepsilon\to 0^+}\tilde{x}_{\varepsilon}=x_0,$$
even if we only have
$\limsup_{\varepsilon\to 0^+} F_{\varepsilon}(\tilde{x}_{\varepsilon})<+\infty$.
\end{teo}

Following again \cite{Ron08} we show the applicability of this abstract result to the
inverse conductivity problem with discontinuous conductivities.

We observe that, in order to guarantee convergence of the regularised solutions to the 
looked for solution, we need to find a metric on the space $X$ such that the following properties are satisfied:
  \begin{enumerate}[1)]
\item the forward operator $\Lambda$ is continuous;
\item $R$ is a regularisation operator for $X$;
\item $\Lambda$ is injective (uniqueness of the inverse problem).
\end{enumerate}

We consider in this subsection $X$ equal to $\mathcal{M}(\lambda_0,\lambda_1)$, or
$\mathcal{M}_{sym}(\lambda_0,\lambda_1)$, or $\mathcal{M}_{scal}(\lambda_0,\lambda_1)$.

On $X$ we consider the metric given by the $L^1$ norm, in all cases. In fact we wish to have a convergence in a rather strong sense, being for instance $H$-convergence
too weak for applications.

Therefore, we take as $Y$ the usual space where
we assume that, for some $p>2$,
in the Dirichlet-to-Neumann case, 
$B_1\subset W^{1-1/p,p}(\partial \Omega)$,
with continuous immersion, and, in the Neumann-to-Dirichlet case, we assume that
$\tilde{B}_1\subset (W^{1-1/p',p'}(\partial \Omega))'_{\ast}$, with continuous immersion.

As a regularisation operator, there are several possibilities. One is to consider a kind of 
total variation regularisation. For instance, 
we define, for any $\sigma\in X$,
$TV(\sigma)$ as the matrix such that
$TV(\sigma)_{ij}=TV(\sigma_{ij})=|D\sigma_{ij}|(\Omega)$ and set
$|\sigma|_{BV(\Omega)}=\| TV(\sigma) \| $ for any $\sigma\in X$.
For any $\sigma\in X$ we define
$$\|\sigma\|_{BV(\Omega)}=\|\sigma\|_{L^1(\Omega)}+|\sigma|_{BV(\Omega)}.$$
Then we may pick as $R$ either $|\cdot|_{BV(\Omega)}$ or $\|\cdot\|_{BV(\Omega)}$.

The total variation regularisation has been widely used in the literature for solving numerically the inverse conductivity problem, for example in \cite{Dob-San94}, with a discretisation method, and in \cite{Chan-Tai, Chu-Chan-Tai}, with level set methods.

Another option is the so-called Mumford-Shah operator. In this case we limit ourselves to scalar conductivities, that is, to $X=\mathcal{M}_{scal}(\lambda_0,\lambda_1)$,
and define, for any $\sigma\in \mathcal{M}_{scal}(\lambda_0,\lambda_1)$,
\begin{equation}\label{MS}
R(\sigma)=\left\{\begin{array}{ll}
\displaystyle{b\int_{\Omega}\|\nabla \sigma\|^2+\mathcal{H}^{N-1}(J(\sigma))} & \text{if }\sigma\in SBV(\Omega)\\
\vphantom{\displaystyle{\int}}+\infty & \text{otherwise}.
\end{array}\right.
\end{equation}
Here $b$ is a positive constant, $\mathcal{H}^{N-1}$ denotes the $(N-1)$-dimensional Hausdorff measure, $J(\sigma)$ is the jump set of $\sigma$, and $SBV$ denotes the space of special functions of bounded variations.
The functional $R$ here defined is referred to as the Mumford-Shah functional and was introduced in the context of image segmentation in \cite{Mum-Sha89}.
We refer, for instance, to \cite{Amb-Fus-Pal} for a detailed discussion on these topics.
The compactness
and semicontinuity theorem for special functions of bounded variation due to Ambrosio, see for instance \cite[Theorem~4.7 and Theorem~4.8]{Amb-Fus-Pal}, guarantees that also in this case $R$ is a regularisation operator for $X$. In the context of inverse problems, and in particular for the inverse conductivity problem, the Mumford-Shah functional has been used as regularisation for the first time in \cite{Ron-San}, with an implementation exploiting the approximation of the Mumford-Shah functional by functionals defined on smoother functions due to Ambrosio and Tortorelli,
\cite{Amb-Tor90,Amb-Tor92}.

We now recall the results in \cite{Ron08}, that immediately follows from the previous abstract results.

\begin{teo}\label{mainapplteo}
Under the previous notation and assumptions, let $\Lambda:X\to Y$ be the forward operator.
Let $R$ be either $|\cdot|_{BV(\Omega)}$ or $\|\cdot\|_{BV(\Omega)}$. If $X=\mathcal{M}_{scal}(\lambda_0,\lambda_1)$, $R$ may be also chosen as in \eqref{MS}.

Let $\sigma_0\in X$ be such that $R(\sigma_0)<+\infty$ and
$\hat{\Lambda}_0=\Lambda(\sigma_0)$. For any $\varepsilon$, $0<\varepsilon\leq \varepsilon_0$, let $\hat{\Lambda}_{\varepsilon}\in Y$ be such that
$\|\hat{\Lambda}_{\varepsilon}-\hat{\Lambda}_0\|\leq \varepsilon$.

Let us fix positive constants $\alpha$, $\gamma$, and $\tilde{a}$, such that
$0<\gamma< \alpha$. For any $\varepsilon$, $0<\varepsilon\leq\varepsilon_0$,
let $F_{\varepsilon}$ be defined as in \eqref{Fdef} and $F_0$ be defined as in \eqref{defF0}.

Then  we have that there exists $\min_X F_{\varepsilon}$, for any $\varepsilon$, $0\leq\varepsilon\leq \varepsilon_0$, and
$$\min_X F_0=\lim_{\varepsilon\to 0^+}\min_X F_{\varepsilon}<+\infty.$$

Let $\{\varepsilon_n\}_{n\in\mathbb{N}}$ be a sequence of positive numbers converging to $0$ as $n\to\infty$.

Let
$\{\tilde{\sigma}_{n}\}_{n\in\mathbb{N}}$ be such that $\limsup_nF_{\varepsilon_n}(\tilde{\sigma}_n)<+\infty$.
Then, up to a subsequence, $\tilde{\sigma}_n$ converges in the $L^1$ norm to $\tilde{\sigma}\in X$
such that $\tilde{\sigma}$ satisfies $\|\Lambda(\tilde{\sigma})-\Lambda(\sigma_0)\|_Y=0$.

Let
$\{\tilde{\sigma}_{n}\}_{n\in\mathbb{N}}$ be such that $\lim F_{\varepsilon_n}(\tilde{\sigma}_n)=
\lim_n\min_X F_{\varepsilon_n}$.
Then, up to a subsequence, $\tilde{\sigma}_n$ converges in the $L^1$ norm to $\tilde{\sigma}\in X$
such that $\tilde{\sigma}$ is a minimizer of $F_0$, that is, in particular, $\|\Lambda(\tilde{\sigma})-\Lambda(\sigma_0)\|_Y=0$
and $R(\tilde{\sigma})=\min\{R(\sigma):\ \sigma\in X\text{ such that }\|\Lambda(\sigma)-\Lambda(\sigma_0)\|_Y=0\}$.
\end{teo}

In dimension $2$ and for scalar conductivities we have the following.

\begin{teo}\label{2dimcase}
Under the notation and assumptions of Theorem~\textnormal{\ref{mainapplteo}},
let us further assume that the space dimension is $2$, that is $N=2$.
We pick $X=\mathcal{M}_{scal}(\lambda_0,\lambda_1)$ and we assume that
either $B_1$ is dense in $H^{1/2}_{\ast}(\partial\Omega)$ or $\tilde{B}_1$ is dense in $H^{-1/2}_{\ast}(\partial\Omega)$, respectively.

Let $\{\tilde{\sigma}_{\varepsilon}\}_{0<\varepsilon\leq\varepsilon_0}$ satisfy $\limsup_{\varepsilon\to 0^+} F_{\varepsilon}(\tilde{\sigma}_{\varepsilon})<+\infty$.
Then we have that
$$\lim_{\varepsilon\to 0^+}\int_{\Omega}|\tilde{\sigma}_{\varepsilon}-\sigma_0|=0.$$
\end{teo}

\bigskip

We notice that, when $N\geq 3$,
even if, recently, a great improvement has been achieved in the uniqueness issue,
still we do not have a uniqueness result for scalar $BV$ or $SBV$ functions. To prove uniqueness, or nonuniqueness, in this case is an extremely interesting and challenging open problem.

We recall that the approach to regularisation developed in \cite{Ron08} has been followed in other works. In \cite{Jia-Ma-Pa} the Mumford-Shah approach has been made slightly more precise, for instance it was proved convergence of the jump sets, and it has been applied to other inverse problems, such as
image deblurring or X-ray tomography.
In \cite{Ji-Ma}, instead, other regularisation strategies for the inverse conductivity problem have been considered, for example the sparsity or 
smoothness penalty was used. In this case the theory for convergence of Tikhonov regularised solutions for nonlinear operators may be used and, in fact, in \cite{Ji-Ma} some convergence estimates were derived.

\section{Numerical approximation and regularisation for the inverse conductivity problem}\label{discretesec}

After the regularisation strategy has been decided, and it has been proved to be effective, the second step is to proceed in finding a suitable numerical approach to solve the regularised minimum problem. For example, in \cite{Ron-San}, the Ambrosio and Tortorelli approximation of the Mumford-Shah functional was used to tackle numerically the minimisation problem. For total variation regularisation, besides the
early paper where a discretisation method, \cite{Dob-San94}, or level set methods, \cite{Chan-Tai,Chu-Chan-Tai}, were used, an interesting analysis of a finite element approximation has been developed in \cite{Ge-Ji-Lu}.

However, the approximations in \cite{Ron-San} and in \cite{Ge-Ji-Lu} have been performed just for the regularised minimum problem, that is, for a fixed regularisation parameter. Instead, we believe that it is very important to study how the approximation parameter (for example the size of the mesh in the finite element approximation) and the regularisation parameter interact. In other words, we wish to find, for a corresponding noise level $\varepsilon$, what are the right regularisation and approximation parameters that allow to prove that the solutions to the approximated regularised minimum problems converge, in a suitable sense, to the looked for solution of the inverse problem. Therefore we wish to include in the convergence analysis developed in \cite{Ron08}, and here recalled in Subsection~\ref{regularsubs}, the approximation of the regularised minimum problem, simultaneously.

Such an approach has been developed for the Ambrosio-Tortorelli approximation of the Mumford-Shah functional in \cite{Che}. For the convenience of the reader we recall the result of \cite{Che} in Subsection~\ref{cherinisubs}.

In the next Subsection~\ref{discretesubsec}, we consider the approximation by finite element discretisation and we investigate how the discretisation parameter should be linked to the noise level and the regularisation parameter. We present here a very simple setting, in future work we will consider a much more general and complete discretisation of the inverse conductivity problem.

Let us begin by introducing the common setting for the whole section.

Throughout this section we fix $\Omega$, a bounded connected open set with Lipschitz boundary, contained in $\mathbb{R}^N$, $N\geq 2$, and two constants $\lambda_0$, $\lambda_1$, with $0<\lambda_0\leq\lambda_1$.

We consider only the case of scalar conductivities, namely we call
$X=\mathcal{M}_{scal}(\lambda_0,\lambda_1)$.

We fix a real number $p>2$.
In the Dirichlet-to-Neumann case, we assume that
$B_1\subset W^{1-1/p,p}_{\ast}(\partial \Omega)$ and $H^{-1/2}_{\ast}(\partial\Omega)\subset B_2$,
with continuous immersions.
In the Neumann-to-Dirichlet case, we assume that
$\tilde{B}_1\subset (W^{1-1/p',p'}(\partial \Omega))'_{\ast}$ and
$H^{1/2}_{\ast}(\partial\Omega)\subset \tilde{B}_2$,
with continuous immersions.

In the Dirichlet-to-Neumann case, we call $Y=\mathcal{L}(B_1,B_2)$ and define
$\Lambda:X\to Y$ as follows
$$\Lambda(\sigma)=DN(\sigma)|_{B_1}:B_1\to B_2.$$

In the Neumann-to-Dirichlet case, we call $Y=\mathcal{L}(\tilde{B}_1,\tilde{B}_2)$ and define
$\Lambda:X\to Y$ as follows
$$\Lambda(\sigma)=ND(\sigma)|_{\tilde{B}_1}:\tilde{B}_1\to \tilde{B}_2.$$

The important fact is the following. We know that $\Lambda:X\to Y$ is H\"older continuous, that is, there exists constant $C_0>0$ and $\beta$, $0<\beta<1$, such that, for any $\sigma$, $\tilde{\sigma}\in X$, we have
\begin{equation}\label{Holderineq}
\|\Lambda(\sigma)-\Lambda(\tilde{\sigma})\|_Y\leq C_0\|\sigma-\tilde{\sigma}\|_{L^1(\Omega)}^{\beta}.
\end{equation}
Here $\beta$ depends on $p$ only and it will play a crucial role in the next analysis.

We consider $\sigma_0\in X$ to be the scalar conductivity in $\Omega$ that we wish to determine and call $\hat{\Lambda}_0=\Lambda(\sigma_0)\in Y$.

Fixed a positive constant $\varepsilon_0$,
for any $\varepsilon$, $0<\varepsilon\leq\varepsilon_0$, let us assume that there exists
$\hat{\Lambda}_{\varepsilon}\in Y$ such that
\begin{equation}\label{error}
\|\hat{\Lambda}_{\varepsilon}-\hat{\Lambda}_0\|_Y\leq\varepsilon.
\end{equation}
Here $\varepsilon$ plays the role of the noise level and $\hat{\Lambda}_{\varepsilon}$ plays the role of the measured Dirichlet-to-Neumann, or Neumann-to-Dirichlet respectively, map.

\subsection{Discrete approximation and regularisation of the inverse conductivity problem}\label{discretesubsec}

Since we wish to consider a discretisation of the problem, we shall make the following assumptions on $\Omega$. We further assume that $\Omega$ is polygonal, that is, $\Omega$ is a polyhedron in $\mathbb{R}^N$.

We use standard conforming piecewise linear finite elements, for which we refer for instance to \cite{Cia}. We shall keep fixed a positive parameter $s$ and a positive constant $h_0$.
We consider, for a fixed parameter $h$, $0<h\leq h_0$, a triangulation $\mathcal{T}_h$ of $\overline{\Omega}$, that is, $\overline{\Omega}=\bigcup_{K\in \mathcal{T}_h}K$, where
each $K\in \mathcal{T}_h$ is a nondegenerate
$N$-simplex, and $\mathcal{T}_h$ satisfies assumption (FEM 1) in \cite[Chapter~2]{Cia}

We then define the finite element space $X_h$
as follows
$$X_h=\{v_h\in C(\overline{\Omega}):\ v_h|_K\in P_1(K)\text{ for any }K\in\mathcal{T}_h\}$$
where $P_1(K)$ is the space of polynomials of order at most $1$ restricted to $K$, that
is, $X_h$ is the finite element space associated to $N$-simplices of type (1). By \cite[Theorem~2.2.3]{Cia} we have that $X_h\subset C(\overline{\Omega})\cap H^1(\Omega)$. It is also clear that $X_{0h}=\{v_h\in X_h:\ v_h|_{\partial\Omega}=0\}$
is contained in $H^1_0(\Omega)$. We call $\Pi_h$ the associated interpolation operator
defined on $C(\overline{\Omega})$.

We assume that $\mathcal{T}_h$ is \emph{regular} in the following classical sense. 
For any $K\in\mathcal{T}_h$ we call $h_K=\mathrm{diam}(K)$ and $\rho_K=\sup
\{\mathrm{diam}(B):\ B\text{ is a ball contained in }K\}$.
Then we assume that
\begin{equation}\label{regularfinite}
h_K\leq h\text{ and }h_K\leq s\rho_K\quad\text{for any }K\in\mathcal{T}_h.
\end{equation}

The following estimate is an immediate consequence of \cite[Theorem~3.1.6]{Cia}.

\begin{teo}\label{Ciarletestimteo}
Let us consider $q\geq 2$ such that $q>N/2$.
Then there exists a constant $C$ such that
for any $u\in W^{2,q}(\Omega)$ we have
\begin{equation}\label{Ciarletest}
\|u-\Pi_h u\|_{W^{1,q}(\Omega)}\leq Ch\|D^2u\|_{L^q(\Omega)}.
\end{equation}
 \end{teo}
 
Our approach to discretisation is the following. As a regularisation operator we consider a total variation penalisation, that is $R$ is given by, for any $\sigma\in X$,
 \begin{equation}\label{Rdefbis}
 R(\sigma)=|\sigma|_{BV(\Omega)}=TV(\sigma)=|D\sigma|(\Omega)\quad\text{or}\quad R(\sigma)=\|\sigma\|_{BV(\Omega)}.
 \end{equation}
Furthermore we shall assume that $R(\sigma_0)<+\infty$, that is, $\sigma_0\in BV(\Omega)$.

For fixed $\tilde{a}$, $0<\gamma<\alpha$,
let us define, for any $\varepsilon$, $0<\varepsilon\leq \varepsilon_0$, and $h$,
$0<h\leq h_0$, the functional $F_{\varepsilon,h}:X\to \mathbb{R}\cup\{+\infty\}$ such that for any $\sigma\in X$
\begin{equation}\label{Fdefbis}
F_{\varepsilon,h}(\sigma)=\left\{\begin{array}{ll}
\displaystyle{\frac{\|\Lambda(\sigma)-\hat{\Lambda}_{\varepsilon}\|^{\alpha}_Y}{\varepsilon^{\gamma}}+\tilde{a}R(\sigma)}&\text{if }
\sigma\in X_h\\
+\infty &\text{otherwise}.
\end{array}
\right.
\end{equation}

Let us immediately notice that any of these functionals $F_{\varepsilon,h}$ admits a minimum over $X$.

We also define $F_0:X\to \mathbb{R}\cup\{+\infty\}$ as before
\begin{equation}\label{defF0bis}
F_0(\sigma)=\left\{\begin{array}{ll}
\tilde{a}R(\sigma) &\text{if }\|\Lambda(\sigma)-\Lambda(\sigma_0)\|_Y=0\\
+\infty &\text{otherwise}
\end{array}
\right.
\end{equation}
 for any $\sigma\in X$.
 
Our aim is to choose $h=h(\varepsilon)$ such that $F_{\varepsilon,h}$ are equicoercive and $\Gamma$-converge, as $\varepsilon\to 0^+$, to $F_0$. 

Therefore, let us consider two sequences $\{\varepsilon_n\}_{n\in\mathbb{N}}\subset (0,\varepsilon_0]$ and $\{h_n\}_{n\in\mathbb{N}}\subset (0,h_0]$ and we assume that
$\lim_n\varepsilon_n=0$. We define $F_n=F_{\varepsilon_n,h_n}$.

The $\Gamma$-liminf inequality is easy to prove. In fact we have the following.

\begin{prop}\label{Gammaliminf}
Let $\{\sigma_n\}_{n\in\mathbb{N}}\subset X$ be such that $\lim_n\sigma_n=\sigma$ in $X$, that is, $\lim_n\|\sigma_n-\sigma\|_{L^1(\Omega)}=0$.

Then
$$F_0(\sigma)\leq\liminf_nF_n(\sigma_n).$$
\end{prop}

\proof{.} If $\liminf_nF_n(\sigma_n)=+\infty$, then there is nothing to prove. We therefore assume, without loss of generality, that
$\liminf_nF_n(\sigma_n)=\lim_nF_n(\sigma_n)<+\infty$. In particular,
for some constant $C$, we have $F_n(\sigma_n)\leq C$ for any $n\in\mathbb{N}$.
Therefore, $\sigma_n\in X_{h_n}$ for any $n\in\mathbb{N}$.

By semicontinuity of the total variation, it is easy to see that
$$\tilde{a}R(\sigma)\leq\liminf_n (\tilde{a}R(\sigma_n))\leq
\liminf_n F_n(\sigma_n).$$

It remains to prove that $\|\Lambda(\sigma)-\Lambda(\sigma_0)\|_Y=0$.
But, by continuity of $\Lambda$, it is easy to see that
\begin{multline*}
\|\Lambda(\sigma)-\Lambda(\sigma_0)\|_Y=\lim_n
\|\Lambda(\sigma_n)-\Lambda(\sigma_0)\|_Y\\\leq
\liminf_n\left[\|\Lambda(\sigma_n)-\hat{\Lambda}_{\varepsilon_n}\|_Y+
\|\hat{\Lambda}_{\varepsilon_n}-\Lambda(\sigma_0)\|_Y\right]\leq
\liminf_n\left[ (C\varepsilon_n^{\gamma})^{1/\alpha}+\varepsilon_n\right]
\end{multline*}
which is obviously equal to $0$.\cvd

\bigskip

The difficult part is to find a recovery sequence. Clearly the existence of the recovery sequence is trivial, by the $\Gamma$-liminf inequality, if $F_0(\sigma)=+\infty$.
Therefore, it is enough to prove the existence of a recovery sequence
when $F_0(\sigma)$ is finite.

\begin{prop}\label{recoverysequenceprop}
We define $h(\varepsilon)=\varepsilon^{3/\beta}$, for any $\varepsilon$, $0<\varepsilon\leq \varepsilon_0$, and recall that $\gamma<\alpha$.

Let $\sigma\in X$ be such that $F_0(\sigma)<+\infty$, that is,
$\sigma\in BV(\Omega)\cap X$ and it satisfies $\|\Lambda(\sigma)-\Lambda(\sigma_0)\|_Y=0$.

Then there exists $\sigma_{\varepsilon}\in X$, for any $\varepsilon$, $0<\varepsilon\leq \varepsilon_0$, such that
$$F_0(\sigma)=\lim_{\varepsilon\to 0^+} F_{\varepsilon,h(\varepsilon)}(\sigma_{\varepsilon}).$$
\end{prop}

Before proving this proposition, let us observe that it implies the following corollary.

\begin{cor}\label{Gamma+coercivecor}
Under the notation and assumptions of Proposition~\textnormal{\ref{recoverysequenceprop}}, we have that
$F_{\varepsilon,h(\varepsilon)}$ $\Gamma$-converges to $F_0$ as $\varepsilon\to 0^+$.

Moreover, the family of functionals $\{F_{\varepsilon,h(\varepsilon)}\}_{0<\varepsilon\leq \varepsilon_0}$ is equicoercive.
\end{cor}

\proof{.} The $\Gamma$-convergence result follows immediately from Propositions~\ref{Gammaliminf} and \ref{recoverysequenceprop}.

About equicoerciveness, we start with the following remark. By Proposition~\ref{recoverysequenceprop}, we can find a constant $C$ such that
\begin{equation}\label{equiminim}
\min_XF_{\varepsilon,h(\varepsilon)}\leq C\quad\text{for any }0<\varepsilon\leq\varepsilon_0.
\end{equation}

Then we define $K=\{\sigma_{\varepsilon}\}_{0<\varepsilon\leq \varepsilon_0}$,
where $\sigma_{\varepsilon}$ is a minimiser for $F_{\varepsilon,h(\varepsilon)}$, for any $0<\varepsilon\leq \varepsilon_0$. We prove that $K$ is relatively compact in $X$.
In fact, by \eqref{equiminim}, we obtain that, for some constant $C_1$,
$R(\sigma_{\varepsilon})\leq C_1$ for any $0<\varepsilon\leq \varepsilon_0$. Then the fact that $K$ is relatively compact follows immediately by the properties of the regularisation operator $R$.\cvd

\bigskip

We now complete the proof of the existence of the recovery sequence.

\proof{ of Proposition~\textnormal{\ref{recoverysequenceprop}}.}
The difficult part is that we need to build the function $\sigma_{\varepsilon}$ in such a way that it belongs to the discrete space $X_{h(\varepsilon)}$, for any $0<\varepsilon\leq\varepsilon_0$.

The construction is the following. First of all we use the fact that $\Omega$ is an extension domain, since it has Lipschitz boundary. Therefore, for any $u\in BV(\Omega)\cap X$, we can find a function $\tilde{u}\in L^{\infty}(\mathbb{R}^N)$ such that
$\tilde{u}|_{\Omega}=u$, $\lambda_0\leq\tilde{u}\leq\lambda_1$ almost everywhere in
$\mathbb{R}^N$,
and, for a constant $C$ depending on $\Omega$ only,
$$|D\tilde{u}|(\mathbb{R}^N)\leq C|Du|(\Omega),$$
and, moreover, $|D\tilde{u}|(\partial\Omega)=0$. This follows immediately by using \cite[Definition~3.20]{Amb-Fus-Pal}, for instance.

We consider our function $\sigma$ and, by a slight abuse of notation, we still call $\sigma$ its extension $\tilde{\sigma}$ to the whole $\mathbb{R}^N$. We fix a positive symmetric mollifier $\eta$,
that is, $\eta\in C^{\infty}_0(B_1(0))$, $\eta\geq 0$, $\int_{B_1(0)}\eta=1$, and such that $\eta(x)$ depends only on $\|x\|$ for any $x\in B_1(0)$. Clearly $\eta\in C^{\infty}_0(\mathbb{R}^N)$ by extending it to $0$ oustide $B_1(0)$.
For any $\delta>0$, we call
$$\eta_{\delta}(x)=\delta^{-N}\eta(x/\delta),\quad x\in\mathbb{R}^N,$$
and, for any $u\in L^1_{loc}(\mathbb{R}^N)$, we call
$$u_{\delta}=\eta_{\delta}\ast u,$$
where as usual $\ast$ denotes the convolution.

We immediately obtain that, for any $\delta>0$,
$\sigma_{\delta}\in C^{\infty}(\mathbb{R}^N)$ and $\lambda_0\leq\sigma_{\delta}\leq\lambda_1$ almost everywhere in $\mathbb{R}^N$. We also have that, locally,
$\sigma_{\delta}$ converges to $\sigma$ as $\delta\to 0^+$ in the $L^1$ norm.
By \cite[Proposition~1.15]{Giu}, we conclude that
\begin{equation}\label{strictconvergence}
\lim_{\delta\to 0^+}\|\sigma_{\delta}-\sigma\|_{L^1(\Omega)}=0\quad\text{and}\quad
\lim_{\delta\to 0^+}|D\sigma_{\delta}|(\Omega)=|D\sigma|(\Omega).
\end{equation}
Actually, by \cite[Lemma~3.24]{Amb-Fus-Pal}, the $L^1$ convergence may be made much more precise. In fact, for a constant $C_1$ depending on $\Omega$ only, we have,
for any $\delta$, $0<\delta\leq 1$,
\begin{equation}\label{firstestima}
\|\sigma-\sigma_{\delta}\|_{L^1(\Omega)}\leq C_1|D\sigma|(\Omega)\delta.
\end{equation}

We choose $q$ as in Theorem~\ref{Ciarletestimteo}.
Since $\sigma_{\delta}\in C^{\infty}(\mathbb{R}^N)$, we obviously have that
$\sigma_{\delta}\in W^{2,q}(\Omega)$, for any $\delta>0$. We need to control its norm in dependence of $\delta$. We notice that, for any multiindex $\alpha$, we have
$D^{\alpha}\sigma_{\delta}=(D^{\alpha}\eta_{\delta})\ast \sigma$. Therefore, for any $\delta$, $0<\delta\leq 1$, and any $p$, $1\leq p\leq +\infty$,
$$\|D^{\alpha}\sigma_{\delta}\|_{L^p(\Omega)}\leq 
C_2\delta^{-|\alpha|},$$
where $C_2$ depends on $\Omega$, $p$, $|\alpha|$, $\eta$, and $\lambda_1$ only.
We conclude that, for a constant $C_3$ depending on $\Omega$, $q$, $\eta$, and $\lambda_1$ only, we have,
for any $0<\delta\leq 1$,
\begin{equation}\label{W2pbound}
\|\sigma_{\delta}\|_{W^{2,q}(\Omega)}\leq C_3\delta^{-2}.
\end{equation}

By Theorem~\ref{Ciarletestimteo}, we obtain that
\begin{equation}\label{Ciarletest2}
\|\sigma_{\delta}-\Pi_h(\sigma_{\delta})\|_{W^{1,q}(\Omega)}\leq C_4h\delta^{-2}
\end{equation}
where $C_4=C_3C$, with $C$ as in \eqref{Ciarletest}. We have that
$\Pi_h(\sigma_{\delta})\in X_h$. Furthermore,
\begin{multline*}
\|\sigma-\Pi_h(\sigma_{\delta})\|_{L^1(\Omega)}\leq
\|\sigma-\sigma_{\delta}\|_{L^1(\Omega)}+
\|\sigma_{\delta}-\Pi_h(\sigma_{\delta})\|_{L^1(\Omega)}\\\leq
C_1|D\sigma|(\Omega)\delta+C_5C_4h\delta^{-2},
\end{multline*}
with $C_5$ depending on $\Omega$ and $q$ only. By picking $\delta=h^{1/3}$, we conclude that, for the constant $C_6=C_1|D\sigma|(\Omega)+C_5C_4$,
\begin{equation}\label{L1estimate}
\|\sigma-\Pi_h(\sigma_{\delta})\|_{L^1(\Omega)}\leq C_6h^{1/3}.
\end{equation}
Furthermore,
\begin{multline*}
|D(\Pi_h(\sigma_{\delta}))|(\Omega)=\int_{\Omega}\|\nabla (\Pi_h(\sigma_{\delta}))\|\\ =\left(
\int_{\Omega}\|\nabla (\Pi_h(\sigma_{\delta}))\|-\int_{\Omega}\|\nabla \sigma_{\delta}\|\right)+\int_{\Omega}\|\nabla \sigma_{\delta}\|.
\end{multline*}
The first term of the right hand side goes to $0$, as $h$, and thus $\delta$, goes to $0$, by \eqref{Ciarletest2}. The second term of the right hand side is exactly
$|D\sigma_{\delta}|(\Omega)$, therefore it goes to $|D\sigma|(\Omega)$, as $h$ goes to $0$, by \eqref{strictconvergence}.

We have therefore constructed, for any $0<h\leq h_0$, $\sigma_h\in X_h$ such that
\begin{equation}\label{finecotruzione}
\|\sigma-\sigma_h\|_{L^1(\Omega)}\leq C_6h^{1/3}\quad\text{and}\quad
\lim_{h\to 0^+}|D\sigma_h|(\Omega)=|D\sigma|(\Omega).
\end{equation}

By \eqref{Holderineq}, we conclude that 
\begin{equation}\label{finecotruzione2}
\|\Lambda(\sigma)-\Lambda(\sigma_h)\|_Y\leq C_0C_6^{\beta}h^{\beta/3}\quad\text{and}\quad
\lim_{h\to 0^+}R(\sigma_h)=R(\sigma).
\end{equation}
Then we easily compute, since $\|\Lambda(\sigma)-\Lambda(\sigma_0)\|_Y=0$,
$$\|\Lambda(\sigma_h)-\hat{\Lambda}_{\varepsilon}\|_Y^{\alpha}\leq \left(\|\Lambda(\sigma_h)-\Lambda(\sigma)\|_Y+
\|\Lambda(\sigma)-\hat{\Lambda}_{\varepsilon}\|_Y\right)^{\alpha}\leq
\left(C_0C_6^{\beta}h^{\beta/3}+\varepsilon\right)^{\alpha}.
$$
If we choose $\gamma<\alpha$ and $h(\varepsilon)$ such that
$h(\varepsilon)=\varepsilon^{3/\beta}$, then we obtain that
$$\lim_{\varepsilon\to 0^+}F_{\varepsilon,h(\varepsilon)}(\sigma_{h(\varepsilon)})=
F_0(\sigma).$$
The proof is concluded.\cvd

\bigskip

By Corollary~\ref{Gamma+coercivecor} and the Fundamental Theorem of $\Gamma$-convergence, Theorem~\ref{fundthm}, the next theorems, which are the main results of the paper, immediately follow.

\begin{teo}\label{mainthm}
Under the previous notation and assumptions, we consider
$X=\mathcal{M}_{scal}(\lambda_0,\lambda_1)$ and let 
 $\Lambda:X\to Y$ be the forward operator.
Let $R$ be either $|\cdot|_{BV(\Omega)}$ or $\|\cdot\|_{BV(\Omega)}$.

Let $\sigma_0\in X$ be such that $R(\sigma_0)<+\infty$ and
$\hat{\Lambda}_0=\Lambda(\sigma_0)$. For any $\varepsilon$, $0<\varepsilon\leq \varepsilon_0$, let $\hat{\Lambda}_{\varepsilon}\in Y$ be such that
$\|\hat{\Lambda}_{\varepsilon}-\hat{\Lambda}_0\|\leq \varepsilon$.

Let us fix positive constants $\alpha$, $\gamma$, and $\tilde{a}$, such that
$0<\gamma< \alpha$. For any $\varepsilon$, $0<\varepsilon\leq\varepsilon_0$,
let $h=h(\varepsilon)$ be given by $h(\varepsilon)=\varepsilon^{3/\beta}$, and
$F_{\varepsilon,h(\varepsilon)}$ be defined as in \eqref{Fdefbis} and $F_0$ be defined as in \eqref{defF0bis}.

Then  we have that there exists $\min_X F_{\varepsilon,h(\varepsilon)}$, for any $\varepsilon$, $0\leq\varepsilon\leq \varepsilon_0$, and
$$\min_X F_0=\lim_{\varepsilon\to 0^+}\min_X F_{\varepsilon,h(\varepsilon)}<+\infty.$$

Let $\{\varepsilon_n\}_{n\in\mathbb{N}}$ be a sequence of positive numbers converging to $0$ as $n\to\infty$.

Let
$\{\tilde{\sigma}_{n}\}_{n\in\mathbb{N}}$ be such that $\limsup_nF_{\varepsilon_n,h(\varepsilon_n)}(\tilde{\sigma}_n)<+\infty$.
Then, up to a subsequence, $\tilde{\sigma}_n$ converges in the $L^1$ norm to $\tilde{\sigma}\in X$
such that $\tilde{\sigma}$ satisfies $\|\Lambda(\tilde{\sigma})-\Lambda(\sigma_0)\|_Y=0$.

Let
$\{\tilde{\sigma}_{n}\}_{n\in\mathbb{N}}$ be such that $\lim_n F_{\varepsilon_n,h(\varepsilon_n)}(\tilde{\sigma}_n)=
\lim_n\min_X F_{\varepsilon_n,h(\varepsilon)}$.
Then, up to a subsequence, $\tilde{\sigma}_n$ converges in the $L^1$ norm to $\tilde{\sigma}\in X$
such that $\tilde{\sigma}$ is a minimiser of $F_0$, that is, in particular, $\|\Lambda(\tilde{\sigma})-\Lambda(\sigma_0)\|_Y=0$
and $R(\tilde{\sigma})=\min\{R(\sigma):\ \sigma\in X\text{ such that }\|\Lambda(\sigma)-\Lambda(\sigma_0)\|_Y=0\}$.
\end{teo}

In the two dimensional case, as before, the result may be made more precise.

\begin{teo}\label{2dimcasebis}
Under the notation and assumptions of Theorem~\textnormal{\ref{mainthm}},
let us further assume that the space dimension is $2$, that is $N=2$.
We assume that
either $B_1$ is dense in $H^{1/2}_{\ast}(\partial\Omega)$ or $\tilde{B}_1$ is dense in $H^{-1/2}_{\ast}(\partial\Omega)$, respectively.

Let $\{\tilde{\sigma}_{\varepsilon}\}_{0<\varepsilon\leq\varepsilon_0}$ satisfy $\limsup_{\varepsilon\to 0^+} F_{\varepsilon,h(\varepsilon)}(\tilde{\sigma}_{\varepsilon})<+\infty$.
Then we have that
$$\lim_{\varepsilon\to 0^+}\int_{\Omega}|\tilde{\sigma}_{\varepsilon}-\sigma_0|=0.$$
\end{teo}

\subsection{Regularisation by the Ambrosio-Tortorelli functionals}\label{cherinisubs}

In this subsection we present the approach to regularisation by using the so-called
Ambrosio-Tortorelli functionals
that was developed in \cite{Che}.
These functionals were introduced in \cite{Amb-Tor90,Amb-Tor92} in order to solve numerically the difficult task of minimising the Mumford-Shah functional. In fact the Ambrosio-Tortorelli functionals are a good approximation, in the $\Gamma$-convergence sense, of the Mumford-Shah functional and they are much easier to compute with.

We recall that $\Omega$ is a fixed bounded connected open set with Lipschitz boundary, contained in $\mathbb{R}^N$, $N\geq 2$.
We consider only the case of scalar conductivities, namely we call
$X=\mathcal{M}_{scal}(\lambda_0,\lambda_1)$,
for two constants $\lambda_0$, $\lambda_1$, with $0<\lambda_0\leq\lambda_1$.

Let us begin with the following definition. We fix a continuous function $V:\mathbb{R}\to\mathbb{R}$ such that $V\geq 0$ everywhere in $\mathbb{R}$ and
$V(t)=0$ if and only if $t=1$. We call $c_V=\int_0^1\sqrt{V(t)}\rmd t$. Let $\psi:\mathbb{R}\to\mathbb{R}$ be a lower semicontinuous, nondecreasing function such that $\psi(0)=0$, $\psi(1)=1$, and $\psi(t)>0$ for any $t>0$. For any $\eta>0$, we fix $o_{\eta}\geq 0$ such that $\lim_{\eta\to 0^+}o_{\eta}/\eta=0$, and we call $\psi_{\eta}=\psi+o_{\eta}$. Given a positive parameter $b$, and for any $\eta>0$, we define the functional
$AT_{\eta}:L^1(\Omega)\times L^1(\Omega)\to [0,+\infty]$ as follows,
for any $(u,v)\in L^1(\Omega)\times L^1(\Omega)$,
\begin{multline}\label{ATdefin}
AT_{\eta}(u,v)\\=\left\{\begin{array}{ll}\displaystyle{\int_{\Omega}\left(b\psi_{\eta}(v)\|\nabla u\|^2+\frac{1}{\eta}V(v)+\eta\|\nabla v\|^2\right)} &\begin{array}{l}\text{if }u\in H^1(\Omega)\cap X\\\text{ and } v\in H^1(\Omega,[0,1])
\end{array}
\\
+\infty &\text{otherwise}.
\end{array}\right.
\end{multline}
Here $H^1(\Omega,[0,1])=\{v\in H^1(\Omega):\ 0\leq v\leq 1\text{ a.e. in }\Omega\}$.

We define a new version of the Mumford-Shah functional as follow. We call $MS:L^1(\Omega)\times L^1(\Omega)\to [0,+\infty]$ the functional such that,
for any $(u,v)\in L^1(\Omega)\times L^1(\Omega)$,
\begin{multline}\label{MSdefin}
MS(u,v)\\=\left\{\begin{array}{ll}\displaystyle{b\int_{\Omega}\|\nabla u\|^2+4c_V\mathcal{H}^{N-1}(J(u))} &\begin{array}{l}\text{if }u\in SBV(\Omega)\cap X\\\text{ and } v=1\text{ a.e. in }\Omega
\end{array}
\\
+\infty &\text{otherwise}.
\end{array}\right.
\end{multline}
Notice that here $v$ just plays the role of a formal variable.

We have the following result.
\begin{teo}\label{AmbTorteo}
We have that, as $\eta\to 0^+$,
$AT_{\eta}$ $\Gamma$-converges to $MS$ in the $L^1(\Omega)\times L^1(\Omega)$ distance.

Moreover, we assume that, for a positive constant $C_0$, we have
$\psi(t)\geq C_0t^2$ for any $t\in [0,1]$.
We consider two sequences $\{\eta_n\}_{n\in\mathbb{N}}\subset (0,1]$, such that
$\lim_n\eta_n=0$, and $\{(u_n,v_n)\}_{n\in\mathbb{N}}\subset L^1(\Omega)\times L^1(\Omega)$.
If there exists a constant $C$ such that
$AT_{\eta_n}(u_n,v_n)\leq C$ for any $n\in\mathbb{N}$, then, as $n\to\infty$,
$v_n$ converges to $v\equiv 1$ in $L^1(\Omega)$ and, up to a subsequence, $u_n$
converges to $u\in X$ in $L^1(\Omega)$.
\end{teo}

\proof{.} The $\Gamma$-convergence follows from \cite{Amb-Tor90,Amb-Tor92}, see also \cite{Bra}.

For the compactness result of the second part, the argument is the following.
The fact that $\lim_nv_n=v\equiv 1$ in $L^1(\Omega)$ is trivial. For the compactness of the sequence $\{u_n\}_{n\in\mathbb{N}}$, first of all we notice that
$\lambda_0\leq u_n\leq \lambda_1$ for any $n\in\mathbb{N}$.
We call $\tilde{V}(t)=\int_0^t\sqrt{V(s)}\rmd s$, for any $t\in [0,1]$. We notice that, for any $t\in [0,1]$, we have
$c_1 t\leq \tilde{V}(t)\leq C_1t$, for some constants $0<c_1<C_1$. Therefore, for any
$t\in (0,1]$, we have
$$\frac{\tilde{V}(t)}{\sqrt{\psi(t)}}\leq C_2.$$

For any $n\in\mathbb{N}$, we define the auxiliary function $w_n=\tilde{V}(v_n)u_n$ and notice that $\|w_n\|_{L^{\infty}(\Omega)}$ is uniformly bounded.
Then $\nabla w_n=\sqrt{V(v_n)}u_n \nabla v_n+\tilde{V}(v_n)\nabla u_n.$
We obtain that
\begin{multline*}\int_{\Omega}\|\nabla w_n\|\leq \|u_n\|_{L^{\infty}(\Omega)}
\left(\int_{\Omega}\frac{1}{\eta_n}V(v_n)\right)^{1/2}\left(\int_{\Omega}\eta_n
\|\nabla v_n\|^2\right)^{1/2}\\
+\left(\int_{\{x\in\Omega:
\ v_n(x)>0\}}\left(\frac{\tilde{V}(v_n)}{\sqrt{\psi(v_n)}})\right)^2\right)^{1/2}\left(\int_{\Omega}
\psi(v_n)
\|\nabla u_n\|^2\right)^{1/2}.
\end{multline*}
We easily conclude that $\{w_n\}_{n\in\mathbb{N}}$ is bounded in
$W^{1,1}(\Omega)$, therefore, up to a subsequence that we do not relabel, we have that
$w_n\to w\in L^1(\Omega)$ and $v_n\to v\equiv 1$, in both cases
in $L^1(\Omega)$ and almost everywhere in $\Omega$.
For almost any $x\in \Omega$, we have that, as $n\to\infty$,
$w_n(x)\to w(x)$ and $v_n(x)\to 1$, thus $\tilde{V}(v_n(x))\to \tilde{V}(1)>0$.
Therefore, for any of these $x\in\Omega$, we have $\lim_nu_n(x)=w(x)/\tilde{V}(1)=u(x)$. We conclude that $u\in X$ and that, up to the same subsequence, as $n\to\infty$,
$u_n$ converges to $u$ almost everywhere in $\Omega$, thus, by the uniform $L^{\infty}$ bound and the Lebesgue theorem, in $L^1(\Omega)$ as well.\cvd

\bigskip

We now consider the following definition.
For fixed $\tilde{a}$, $0<\gamma<\alpha$,
let us define, for any $\varepsilon$, $0<\varepsilon\leq \varepsilon_0$, and $\eta$,
$0<\eta\leq \eta_0$, the functional $F_{\varepsilon,\eta}:X\times L^1(\Omega)\to \mathbb{R}\cup\{+\infty\}$ such that, for any $(\sigma,v)\in X\times L^1(\Omega)$, we have
\begin{equation}\label{Fdefter}
F_{\varepsilon,\eta}(\sigma,v)=
\displaystyle{\frac{\|\Lambda(\sigma)-\hat{\Lambda}_{\varepsilon}\|^{\alpha}_Y}{\varepsilon^{\gamma}}+\tilde{a}AT_{\eta}(\sigma,v)}.
\end{equation}

We also define $F_0:X\times L^1(\Omega)\to \mathbb{R}\cup\{+\infty\}$ as follows, 
for any $(\sigma,v)\in X\times L^1(\Omega)$,
\begin{equation}\label{defF0ter}
F_0(\sigma,v)=\left\{\begin{array}{ll}
\tilde{a}MS(\sigma,v) &\text{if }\|\Lambda(\sigma)-\Lambda(\sigma_0)\|_Y=0\\
+\infty &\text{otherwise}
\end{array}
\right.
\end{equation}
where $MS$ is defined in \eqref{MSdefin}.
We notice that, equivalently, we can consider $\tilde{F}_0:X\to \mathbb{R}\cup\{+\infty\}$ such that,
for any $\sigma\in X$,
\begin{multline}\label{defF0quater}
\tilde{F}_0(\sigma)\\=\left\{\begin{array}{ll}
\displaystyle{\tilde{a}\left(b\int_{\Omega}\|\nabla \sigma\|^2+4c_V\mathcal{H}^{N-1}(J(\sigma))\right)} &\begin{array}{l}\text{if }\sigma\in SBV(\Omega)\cap X\\\text{ and }
\|\Lambda(\sigma)-\Lambda(\sigma_0)\|_Y=0
\end{array}
\\
+\infty &\text{otherwise}.
\end{array}
\right.
\end{multline}

\begin{oss}
We notice that $F_0$, or equivalently $\tilde{F}_0$, admits a minimum over $X\times L^1(\Omega)$, or $X$ respectively. Notice that
$(\tilde{\sigma},\tilde{v})$ is a minimiser for $F_0$ if and only if
$\tilde{\sigma}$ is a minimiser for $\tilde{F}_0$ and $\tilde{v}\equiv 1$.
Moreover, 
any of these functionals $F_{\varepsilon,\eta}$ admits a minimum over $X\times L^1(\Omega)$ provided $o_{\eta}>0$.
\end{oss}

We shall need the following definition.

\begin{defin}\label{Minko+adm}
For any Borel set $E\subset \mathbb{R}^N$, we define its $(N-1)$-\emph{dimensional Minkowski content} as
$$\mathcal{M}^{N-1}(E)=\lim_{\delta\to 0^+}\frac{1}{2\delta}|\{x\in\mathbb{R}^N:\ \mathrm{dist}(x,E)<\delta\}|,$$
provided the limit exists.

We say that a conductivity $\sigma\in X$ is \emph{admissible} if $\sigma\in SBV(\Omega)$, and it satisfies
$$\int_{\Omega}\|\nabla \sigma\|^2+\mathcal{H}^{N-1}(J(\sigma))<+\infty\quad\text{and}\quad\mathcal{M}^{N-1}(J(\sigma))=\mathcal{H}^{N-1}(J(\sigma)).$$
\end{defin}

With this definition at hand, we consider the following lemma.

\begin{lem}\label{convlemma}
Let $\sigma\in X$ be admissible in the sense of Definition~\textnormal{\ref{Minko+adm}}. Then we can find $(\sigma_{\eta},v_{\eta})\in L^1(\Omega)\times L^1(\Omega)$, $0<\eta\leq \eta_0$, such that, for some constant $C$,
\begin{equation}\label{convrecovery}
\|\sigma_{\eta}-\sigma\|_{L^1(\Omega)}\leq C\eta \quad\text{and}\quad\lim_{\eta\to 0^+}
\|v_{\eta}-1\|_{L^1(\Omega)}=0,
\end{equation}
and
\begin{equation}\label{recoveryCher}
\lim_{\eta\to 0^+}AT_{\eta}(\sigma_{\eta},v_{\eta})=MS(\sigma,1).
\end{equation}
\end{lem}

\proof{.} Let $\phi:\mathbb{R}\to\mathbb{R}$ be a $C^{\infty}$ function that is nondecreasing and such that $\phi(t)=0$ for any $t\leq 1/8$ and
$\phi(t)=1$ for any $t\geq 7/8$.

For the time being, we consider the case in which $o_{\eta}>0$ and we define
$\xi_{\eta}=\sqrt{\eta o_\eta}$.

We define, for any $\eta$, $0<\eta\leq\eta_0$, and any $x\in\Omega$, 
$$\phi_{\eta}(x)=\phi\left(\frac{\mathrm{dist}(x,J(\sigma))}{\xi_{\eta}}\right).$$ 
Then we define
$$\sigma_{\eta}=\phi_{\eta}\sigma+(1-\phi_{\eta})\lambda_0$$
and, for any $x\in\Omega$, and any $\delta>0$,
$$v_{\eta}^{\delta}(x)=\left\{
\begin{array}{ll}
0 & \text{if }\mathrm{dist}(x,J(\sigma))< \xi_{\eta}\\
\displaystyle{v^{\delta}\left(\frac{\mathrm{dist}(x,J(\sigma))-\xi_{\eta}}{\eta}\right)}& \text{if }\xi_{\eta}\leq\mathrm{dist}(x,J(\sigma))< \xi_{\eta}+T\eta\\
1 & \text{if }\mathrm{dist}(x,J(\sigma))\geq \xi_{\eta}+T\eta.
\end{array}
\right.$$
Here the function $v=v^{\delta}$, and the constant $T>0$, are chosen in such a way $v\in C^1([0,T])$, $v(0)=0$, $v(T)=1$, and
$$\int_{0}^T(V(v)+|v'|^2) \leq 2c_V+\delta.$$

We call, for any positive $r>0$, $S_r=\{x\in\mathbb{R}^N:\ \mathrm{dist}(x,J(\sigma))<r\}$.
First of all we notice that, for some constant $C$,
$$\|\sigma_{\eta}-\sigma\|_{L^1(\Omega)}\leq (\lambda_1-\lambda_0)|S_{\xi_{\eta}}\cap \Omega|\leq C\xi_{\eta}\leq C\eta\quad\text{for any }0<\eta\leq\eta_0.$$

Since $|S_{\xi_{\eta}+T\eta}|\to 0$, as $\eta\to 0^+$, we also deduce that
$v_{\eta}^{\delta}\to 1$ almost everywhere in $\Omega$ and in $L^1(\Omega)$ as well, in a completely independent way from the constant $\delta$ that may be chosen as depending from $\eta$.
 
Then we can compute, since obviously we have that $\sigma_{\eta}\in H^1(\Omega)\cap X$ and $v_{\eta}\in H^1(\Omega,[0,1])$,
\begin{multline*}
AT_{\eta}(\sigma_{\eta},v_{\eta}^{\delta})=\int_{\Omega}\left(b\psi_{\eta}(v_{\eta}^{\delta})\|\nabla \sigma_{\eta}\|^2+\frac{1}{\eta}V(v_{\eta}^{\delta})+\eta\|\nabla v_{\eta}^{\delta}\|^2\right)\\=
b\int_{\Omega\backslash S_{\xi_{\eta}}}\psi_{\eta}(v_{\eta}^{\delta})\|\nabla\sigma\|^2+
bo_{\eta}\int_{S_{\xi_{\eta}}}\|\nabla\sigma_{\eta}\|^2+\frac{1}{\eta}V(0)|S_{\xi_{\eta}}\cap\Omega|\\+
\int_{(S_{\xi_{\eta}+T\eta}\backslash S_{\xi_{\eta}})\cap\Omega}\left(\frac{1}{\eta}V(v_{\eta}^{\delta})+\eta\|\nabla v_{\eta}^{\delta}\|^2\right).
\end{multline*}
Since $v_{\eta}^{\delta}$ converges to $1$ almost everywhere in $\Omega$, it is straightforward to see that the first three terms converge, as $\eta\to 0^+$,
to $\int_{\Omega}b\|\nabla \sigma\|^2$, in a completely independent way from the constant $\delta$ that may be chosen as depending from $\eta$.

By the coarea formula, the definition of the Minkowski content, and the properties of 
$\sigma$, we can prove that
\begin{multline*}
\lim_{\eta\to 0^+}\int_{S_{\xi_{\eta}+T\eta}\backslash S_{\xi_{\eta}}}\left(\frac{1}{\eta}V(v_{\eta}^{\delta})+\eta\|\nabla v_{\eta}^{\delta} \|^2\right)\\=2\left(\int_{0}^TV(v^{\delta})+|(v^{\delta})'|^2\right)\mathcal{M}^{N-1}(J(\sigma)).
\end{multline*}
Since 
$\mathcal{M}^{N-1}(J(\sigma)) =\mathcal{H}^{N-1}(J(\sigma))$, we easily deduce that, even if $o_{\eta}=0$,
$$\limsup_{\eta\to 0^+}AT_{\eta}(\sigma_{\eta},v_{\eta}^{\delta})\leq \int_{\Omega}b\|\nabla \sigma\|^2+(4c_V+2\delta)\mathcal{H}^{N-1}(J(\sigma)).$$
It is then easy to choose $\delta=\delta(\eta)$ and define $v_{\eta}=v_{\eta}^{\delta(\eta)}$, $0<\eta\leq\eta_0$, in such a way that
$$\limsup_{\eta\to 0^+}AT_{\eta}(\sigma_{\eta},v_{\eta})\leq \int_{\Omega}b\|\nabla \sigma\|^2+4c_V\mathcal{H}^{N-1}(J(\sigma)).$$

Clearly $\lim_{\eta\to 0^+}\|v_{\eta}-1\|_{L^1(\Omega)}=0$. Hence,
by the corresponding $\Gamma$-liminf inequality proved in \cite[Proposition~4.5]{Bra},
we conclude that
$$\lim_{\eta\to 0^+}AT_{\eta}(\sigma_{\eta},v_{\eta})= \int_{\Omega}b\|\nabla \sigma\|^2+4c_V\mathcal{H}^{N-1}(J(\sigma))=MS(\sigma,1).$$
Thus the proof is complete.\cvd

\bigskip

We are ready to state the final convergence result.

\begin{teo}\label{Cherinithm}
Under the previous assumptions, let us assume that $\sigma_0$ is admissible in the sense of Definition~\textnormal{\ref{Minko+adm}}.  Let us also assume that, for a positive constant $C_0$, we have
$\psi(t)\geq C_0t^2$ for any $t\in [0,1]$.

If we pick $\eta=\eta(\varepsilon)=\varepsilon^{1/\beta}$, and we call
$F_{\varepsilon}=F_{\varepsilon,\eta(\varepsilon)}$ as in \eqref{Fdefter}, then we obtain that
$$\min_{X}\tilde{F}_0\leq \liminf_{\varepsilon\to 0^+}\left(\inf_{X\times L^1(\Omega)}F_{\varepsilon}\right)\leq \limsup_{\varepsilon\to 0^+}\left(\inf_{X\times L^1(\Omega)}F_{\varepsilon}\right)<+\infty.$$

Furthermore, let us consider two sequences $\{\varepsilon_n\}_{n\in\mathbb{N}}\subset (0,\varepsilon_0]$, such that
$\lim_n\varepsilon_n=0$, and $\{(\sigma_n,v_n)\}_{n\in\mathbb{N}}\subset X\times L^1(\Omega)$.

If there exists a constant $C$ such that
$F_{\varepsilon_n}(\sigma_n,v_n)\leq C$ for any $n\in\mathbb{N}$, then, as $n\to\infty$,
$v_n$ converges to $v\equiv 1$ in $L^1(\Omega)$ and, up to a subsequence, $\sigma_n$
converges to $\tilde{\sigma}\in X$ in $L^1(\Omega)$. Moreover,
$\tilde{\sigma}\in SBV(\Omega)$, $MS(\tilde{\sigma},1)$ is finite, and
$\|\Lambda(\tilde{\sigma})-\Lambda(\sigma_0)\|_Y=0$.
Finally, if $N=2$ and we assume that
either $B_1$ is dense in $H^{1/2}_{\ast}(\partial\Omega)$ or $\tilde{B}_1$ is dense in $H^{-1/2}_{\ast}(\partial\Omega)$, respectively, the whole sequence $\sigma_n$ converges, as $n\to\infty$, to $\sigma_0$ in $L^1(\Omega)$.
\end{teo}

\proof{.} First of all, by applying Lemma~\ref{convlemma} to $\sigma=\sigma_0$, we conclude that 
$$\limsup_{\varepsilon\to 0^+}\left(\inf_{X\times L^1(\Omega)}F_{\varepsilon}\right)<+\infty.$$
In fact, for any $0<\varepsilon\leq\varepsilon_0$, we have
$$\|\Lambda(\sigma_{\eta(\varepsilon)})-\hat{\Lambda}_{\varepsilon}\|_Y\leq
\|\Lambda(\sigma_{\eta(\varepsilon)})-\Lambda(\sigma_0)\|_Y+\|
\Lambda(\sigma_0)-
\hat{\Lambda}_{\varepsilon}\|_Y\leq C(\eta(\varepsilon))^{\beta}+\varepsilon\leq C_1\varepsilon,$$
for some constants $C$ and $C_1$.

By the
$\Gamma$-limif inequality, \cite[Proposition~4.5]{Bra}, and the compactness stated in the second part of Theorem~\ref{AmbTorteo},
we can immediately prove that
$$\min_{X}\tilde{F}_0\leq \liminf_{\varepsilon\to 0^+}\left(\inf_{X\times L^1(\Omega)}F_{\varepsilon}\right).$$

The second part of the theorem follows immediately, again by exploiting the compactness result in Theorem~\ref{AmbTorteo}.\cvd


\begin{thebibliography}{99}


\bibitem{Ale-Cab}
G.~Alessandrini and E.~Cabib,
\emph{EIT and the average conductivity},
J. Inverse Ill-Posed Probl. \textbf{16} (2008) 727--736.

\bibitem{A-dH-G-Sin}
G.~Alessandrini, M.~V.~de Hoop, R.~Gaburro and E.~Sincich, 
\emph{Lipschitz stability for the electrostatic inverse boundary value problem with piecewise linear conductivities}, preprint (2015).

\bibitem{Ale-Ves}
G.~Alessandrini and S.~Vessella,
\emph{Lipschitz stability for the inverse conductivity problem},
Adv. in Appl. Math. \textbf{35} (2005) 207--241. 

\bibitem{All}
G.~Allaire,
\emph{Shape Optimization by the Homogenization Method},
Springer-Verlag, New York, 2002.

\bibitem{Amb-Fus-Pal}
L.~Ambrosio, N.~Fusco and D.~Pallara,
\emph{Functions of Bounded Variation and Free Discontinuity Problems},
Clarendon Press, Oxford, 2000.

\bibitem{Amb-Tor90}
L.~Ambrosio and V.~M.~Tortorelli,
\emph{Approximation of functionals depending on jumps by elliptic
functionals via $\Gamma$-convergence},
Comm. Pure Appl. Math. \textbf{43} (1990) 999--1036.

\bibitem{Amb-Tor92}
L.~Ambrosio and V.~M.~Tortorelli,
\emph{On the approximation of free discontinuity problem},
Boll. Un. Mat. Ital. \textbf{6-B} (1992) 105--123.

\bibitem{Ast-Pai}
K.~Astala and L.~P\"aiv\"arinta,
\emph{Calder\'on's inverse conductivity problem in the plane},
Ann. of Math. (2) \textbf{163} (2006) 265--299.

\bibitem{Ast-Pai-Las}
K.~Astala, L.~P\"aiv\"arinta and M.~Lassas,
\emph{Calder\'on's inverse problem for anisotropic conductivity in the plane},
Comm. Partial Differential Equations \textbf{30} (2005) 207--224.

\bibitem{Bra}
A.~Braides,
\emph{Approximation of Free-Discontinuity Problems},
Springer-Ver\-lag, Berlin, 1998.

\bibitem{Cal}
A.~P.~Calder\'on,
\emph{On an inverse boundary value problem},
in \emph{Seminar on numerical analysis and its applications to continuum
physics,}
Sociedade Brasileira de Matem\'atica, Rio de Janeiro, 1980, pp.~65--73.

\bibitem{Car-Rog}
P.~Caro and K.~M.~Rogers,
\emph{Global uniqueness for the Calder\'on's problem with Lipschitz conductivities},
Forum Math. Pi \textbf{4} (2016) e2 (28 pp). 

\bibitem{Chan-Tai}
T.~F.~Chan and X.-C.~Tai,
\emph{Level set and total variation regularization for elliptic inverse problems with discontinuous coefficients},
J. Comput. Phys. \textbf{193} (2003) 40--66.

\bibitem{Che}
A.~Cherini,
\emph{Regolarizzazione del problema inverso della conduttivit\`a per conduttivit\`a discontinue},
Master Thesis, Universit\`a degli Studi di Trieste, academic year 2008/2009.

\bibitem{Chu-Chan-Tai}
E.~T.~Chung, T.~F.~Chan and X.-C.~Tai,
\emph{Electrical impedance tomography using level set representation and total variational regularization},
J. Comput. Phys. \textbf{205} (2005) 357--372.

\bibitem{Cia}
P.~G.~Ciarlet,
\emph{The Finite Element Method for Elliptic Problems},
North-Holland, Amsterdam, 1978.

\bibitem{DM}
G.~Dal Maso,
\emph{An Introduction to $\Gamma$-convergence},
Birkh\"auser, Boston, 1993.

\bibitem{Des}
A.~Dessenibus,
\emph{Continuit\`a della mappa
Dirichlet-a-Neumann rispetto alla
conduttivit\`a},
Master Thesis, Universit\`a degli Studi di Trieste, academic year 2014/2015.

\bibitem{Dob-San94}
D.~C.~Dobson and F.~Santosa,
\emph{An image-enhancement technique for electrical impedance tomography},
Inverse Problems \textbf{10} (1994) 317--334.

\bibitem{Eng-et-al}
H.~W.~Engl, M.~Hanke and A.~Neubauer,
\emph{Regularization of Inverse Problems},
Kluwer Academic Publishers, Dordrecht Boston London, 1996.

\bibitem{Eng-Kun-Neu}
H.~W.~Engl, K.~Kunisch and A.~Neubauer,
\emph{Convergence rates for Tikhonov regularisation of nonlinear ill-posed problems},
Inverse Problems \textbf{5} (1989) 523--540. 

\bibitem{Far-Kur-Rui}
D.~Faraco, Y.~Kurylev and A.~Ruiz,
$G$-\emph{convergence}, \emph{Dirichlet to Neumann maps and invisibility},
J. Funct. Anal. \textbf{267} (2014) 2478--2506.

\bibitem{Ge-Ji-Lu}
M.~Gehre, B.~Jin and X.~Lu,
\emph{An analysis of finite element approximation in electrical impedance tomography},
Inverse Problems \textbf{30} (2014) 045013 (24 pp). 

\bibitem{Giu}
E.~Giusti,
\emph{Minimal Surfaces and Functions of Bounded Variation},
Bir\-kh\"auser, Basel, 1984.

\bibitem{Hab}
B.~Haberman,
\emph{Uniqueness in Calder\'on's problem for conductivities with unbounded gradient}, 
Comm. Math. Phys. \textbf{340}  (2015) 639--659.

\bibitem{Hab-Tat}
B.~Haberman and D.~Tataru,
\emph{Uniqueness in Calder\'on's problem with Lipschitz conductivities},
Duke Math. J. \textbf{162} (2013) 496--516.

\bibitem{Jia-Ma-Pa}
M.~Jiang, P.~Maass and T.~Page,
\emph{Regularizing properties of the Mumford-Shah functional for imaging applications},
Inverse Problems \textbf{30} (2014) 035007 (17 pp).

\bibitem{Ji-Ma}
B.~Jin and P.~Maass, \emph{ An analysis of electrical impedance tomography with applications to Tikhonov regularization},
ESAIM Control Optim. Calc. Var. \textbf{18} (2012) 1027--1048.

\bibitem{Koh-Vog84:1}
R.~Kohn and M.~Vogelius,
\emph{Determining conductivity
by boundary measurements},
Comm. Pure Appl. Math. \textbf{37} (1984) 289--298.

\bibitem{Koh-Vog85}
R.~V.~Kohn and M.~Vogelius,
\emph{Determining conductivity
by boundary measurements} II. \emph{Interior results},
Comm. Pure Appl. Math. \textbf{38} (1985) 643--667.

\bibitem{LMP}
M.~Lukaschewitsch, P.~Maass and M.~Pidcock,
\emph{Tikhonov regularization for electrical impedance tomography on unbounded domains},
Inverse Problems \textbf{19} (2003) 585--610. 


\bibitem{Mum-Sha89}
D.~Mumford and J.~Shah,
\emph{Optimal approximations by piecewise smooth functions and associated
variational problems},
Comm. Pure Appl. Math. \textbf{42} (1989) 577--685.

\bibitem{Mur-Tar2}
F.~Murat and L.~Tartar,
\emph{Calculus of Variations and Homogenization},
in A.~Cherkaev and R.~Kohn eds., \emph{Topics in the Mathematical Modelling of Composite Materials},
Birkh\"auser, Boston, 1997, pp.~139--173.

\bibitem{Mur-Tar1}
F.~Murat and L.~Tartar,
\emph{$H$-convergence},
in A.~Cherkaev and R.~Kohn eds., \emph{Topics in the Mathematical Modelling of Composite Materials},
Birkh\"auser, Boston, 1997, pp.~21--43.

\bibitem{Nac}
A.~I.~Nachman,
\emph{Global uniqueness for a two-dimensional inverse boundary value
problem},
Ann. of Math. (2) \textbf{143} (1996) 71--96.

\bibitem{Riv-Bar-Ob}
C.~Rivas, P.~Barbone and A.~Oberai,
\emph{Divergence of finite element formulations for inverse problems treated as optimization problems},
 J. Phys.: Conf. Ser. \textbf{135} (6th International Conference on Inverse Problems in Engineering: Theory and Practice) (2008) 012088 (8 pp).

\bibitem{Ron06}
L.~Rondi, \emph{A remark on a paper by Alessandrini and Vessella}, Adv. in Appl. Math. \textbf{36} (2006) 67--69.

\bibitem{Ron08}
L.~Rondi,
\emph{On the regularization of the inverse conductivity problem with discontinuous conductivities}, Inverse Probl. Imaging \textbf{2} (2008) 397--409.

\bibitem{Ron15}
L.~Rondi, \emph{Continuity properties of Neumann-to-Dirichlet maps with respect to the $H$-convergence of the coefficient matrices}, Inverse Problems \textbf{31} (2015) 045002 (24 pp). 

\bibitem{Ron-San}
L.~Rondi and F.~Santosa,
\emph{Enhanced electrical impedance
tomography \emph{via} the Mumford-Shah functional},
ESAIM Control Optim. Calc. Var. \textbf{6} (2001) 517--538.

\bibitem{Syl}
J.~Sylvester,
\emph{An anisotropic inverse boundary value problem},
Comm. Pure Appl. Math. \textbf{43} (1990) 201--232. 

\bibitem{Syl-Uhl87}
J.~Sylvester and G.~Uhlmann,
\emph{A global uniqueness theorem for an inverse boundary value problem},
Ann. of Math. (2) \textbf{125} (1987) 153--169.

\end{thebibliography}
\end{document}